\documentclass[11pt,british,reqno]{amsart}

\usepackage[T1]{fontenc}

\usepackage{babel,eucal,url,amssymb,stmaryrd,booktabs,enumerate,amscd,paralist,cite}

\theoremstyle{plain}
\newtheorem{lemma}{Lemma}[section]
\newtheorem{definition}[lemma]{Definition}
\newtheorem{proposition}[lemma]{Proposition}

\newtheorem{theorem}[lemma]{Theorem}

\newtheorem{example}[lemma]{Example}

\newtheorem*{ack}{Acknowledgements}

\newcommand{\Lie}[1]{\operatorname{\textsl{#1}}}
\newcommand{\lie}[1]{\operatorname{\mathfrak{#1}}}
\newcommand{\GL}{\Lie{GL}}

\newcommand{\SO}{\Lie{SO}}
\newcommand{\so}{\lie{so}}

\newcommand{\Spin}{\Lie{Spin}}
\newcommand{\SU}{\Lie{SU}}

\newcommand{\Un}{\Lie{U}}
\newcommand{\un}{\lie{u}}
\newcommand{\Gtwo}{\ifmmode{{\rm G}_2}\else{${\rm G}_2$}\fi}

 \newcommand{\cyclic}{\mathop{\kern0.9ex{{+}\kern-2.2ex\raise-.28ex\hbox{\Large\hbox
 {$\circlearrowright$}}}}}

\newcommand{\real}[1]{\left\llbracket #1 \right\rrbracket}

\newcommand{\Wc}[1]{\mathcal W_{#1}}
\newcommand{\End}{\mathop{\mbox{\rm End}}}

\def\sideremark#1{\ifvmode\leavevmode\fi\vadjust{\vbox to0pt{\vss
 \hbox to 0pt{\hskip\hsize\hskip1em
 \vbox{\hsize2.5cm\tiny\raggedright\pretolerance10000
 \noindent #1\hfill}\hss}\vbox to8pt{\vfil}\vss}}}%

\newfont{\eusm}{eusm10 scaled \magstep1}
\newfont{\eusmiii}{eusm10 scaled \magstep3}

\newcommand{\comp}{\makebox[7pt]{\raisebox{1.5pt}{\tiny $\circ$}}}
\newcommand{\RR}{\mbox{{\sl I}}\!\mbox{{\sl R}}}

\newcommand{\ve}{^{vert}}
\newcommand{\h}{^{hor}}

\title{\bf Harmonicity of sections of sphere bundles}

\author{J.~C.~Gonz{\'a}lez-D{\'a}vila}
\address[J.~C.~Gonz{\'a}lez-D{\'a}vila]{Department of Fundamental Mathematics\\
  University of La Laguna\\ 38200 La Laguna, Tenerife, Spain}
\email{jcgonza@ull.es}

\author{F.~Mart\'\i n~Cabrera}
\address[F.~Mart\'\i n~Cabrera]{Department of Fundamental Mathematics\\
  University of La Laguna\\ 38200 La Laguna, Tenerife, Spain}
\email{fmartin@ull.es}

\author{M.~Salvai}
\address[M.~Salvai]{famaf-ciem\\
  Ciudad Universitaria\\ 5000 C{\'o}rdoba, Argentina}
\email{salvai@mate.uncor.edu}

\date{\today}

\begin{document}
\maketitle

\begin{abstract}{\indent}
We consider the energy functional on the space of sections of a
sphere  bundle over a Riemannian manifold $(M,\langle \cdot, \cdot
\rangle)$ equipped with the Sasaki metric and we discuss  the
characterising condition for critical points. Likewise, we provide
a useful method for computing the tension field in some particular
situations. Such a method is shown to be adequate for many  tensor
fields defined on manifolds $M$ equipped with a $G$-structure
compatible with $\langle \cdot, \cdot \rangle$. This leads to the
construction of a lot of new examples of differential forms which
are harmonic sections or determine a harmonic map from $(M,\langle
\cdot, \cdot \rangle)$ into its sphere  bundle. \vspace{4mm}

 \vspace{2mm}

\end{abstract}

 \tableofcontents

\section{Introduction}{\indent}
The energy of a map between Riemannian manifolds is a functional
which  has been widely studied by diverse authors
\cite{EeLe1,EeLe2,Ur}. Critical points for the energy functional
are called \emph{harmonic maps} and have been characterised by
Eells and Sampson \cite{EeSa} as maps whose {\it tension field}
vanishes.

For a Riemannian manifold $(M,\langle \cdot , \cdot \rangle)$, we
denote by $(T_1 M, \langle \cdot , \cdot \rangle^S)$ its unit
tangent bundle equipped with the Sasaki metric $\langle \cdot ,
\cdot \rangle^S$ \cite{Bl}. Thinking of unit vector fields as
sections $M \to T_1 M$, if $M$ is  compact and oriented,  one can
consider the energy functional as defined on the set $\mathfrak
X_1(M)$ of unit vector fields. Critical points for this functional
give rise to the notion of \emph{harmonic unit vector field}. The
condition characterising harmonic unit vector fields has been
obtained by Wiegmink \cite{Wie1} (see also Wood's paper
\cite{Wood}). It is interesting to note that harmonic unit vector
fields are not necessarily critical points for the energy
functional on the space of all maps $M \to T_1 M$. A harmonic unit
vector field will also be harmonic map if and only if certain
condition involving the curvature of $M$ is satisfied.

In this paper we consider the energy functional defined on the
space of sections of Riemannian vector bundles $\mathbb E \to M$,
equipped with a metric which generalizes  the Sasaki metric. In
\cite{GGV}, Gil-Medrano et al. considered the energy  functional
defined on $(r,s)$-tensorial bundles on $M$ which is a particular
case of the vector bundles considered here. The characterising
condition of critical points for the energy functional on the
space of sections of  sphere  bundles was shown in \cite{S}. This
gives rise to the notion of \emph{harmonic section of a sphere
bundle}. Additionally,  we  analyse when such harmonic sections
are also harmonic maps. In particular, we will show a method,
mainly based on Lemma \ref{compequi}, for computing the tension
field which is adequate for many situations. Concretely, when we
consider Riemannian manifolds $(M,\langle \cdot , \cdot \rangle)$
equipped with some $G$-structure compatible with the metric
$\langle \cdot , \cdot \rangle$. It is well known that, associated
with many such $G$-structures, there are tensors $\Psi$ of
constant length stabilised by the action of the Lie group $G$. For
many of those $\Psi$, our method proves to be efficient for
computing their corresponding tension fields.

Some of the tensors playing the r\^{o}le of $\Psi$ are: the
K\"{a}hler form of an almost Hermitian structure, the fundamental
three-form of a $G_2$-structure, the real  and imaginary parts of
the complex volume form of a special almost Hermitian structure,
etc. So that we analyse the harmonicity of sections $\Psi$ of
sphere  bundles both as sections and as maps. If the intrinsic
torsion $\xi^G$ of the $G$-structure vanishes, the harmonicity in
both senses trivially follows. Therefore, we study examples
defined on manifolds equipped with a $G$-structure such that
$\xi^G \neq 0$, but with a geometry strongly conditioned by the
$G$-structure. Thus, we have found examples of harmonic sections
and harmonic maps into sphere   bundles defined on manifolds with
$G$-structure  such that their respective intrinsic torsions  have
to be contained in a one-dimensional $G$-module (in many cases,
this implies that the manifold is Einstein): nearly K\"{a}hler
$6$-manifolds, $7$-manifolds with nearly parallel $G_2$-structure,
Sasakian manifolds, Kenmotsu manifolds, etc.

Finally, we focus  attention on manifolds equipped with a locally
conformal parallel $G$-structure. The geometry of such manifolds
is very conditioned by a closed one-form, called the {\it Lee
form}. Thus, for such geometries, we have found tensor fields
which are harmonic sections of sphere   bundles. Furthermore, for
some of them, if the Lee form is parallel, then they are also
harmonic maps.

\begin{ack} {\rm  The first author is supported by a grant from
MEC (Spain), project MTM2007-65852, the second one by a grant from
MEC (Spain), project MTM2007-66375 and the third one by Conicet,
Secyt-UNC and Foncyt.}
\end{ack}

 \section{Preliminaries}{\indent}\setcounter{equation}{0}
The {\em energy} of a map $f:(M, \langle \cdot , \cdot \rangle_M
)\to (N, \langle \cdot , \cdot \rangle_{N})$ between Riemannian
manifolds, $M$ being compact and oriented, is the integral
\begin{equation}\label{f1}
 {\mathcal E}(f)=\frac{\textstyle 1}{\textstyle2}\int_{M}\|f_{*}\|^{2}dv,
\end{equation}
where $\|f_{*}\|$ is the norm of the differential $f_{*}$ of $f$
with respect to the metrics $ \langle \cdot , \cdot \rangle_M$, $
\langle \cdot , \cdot \rangle_N$ and $dv$ denotes the volume form
on $(M,\langle \cdot , \cdot \rangle_M)$.  On the domain of a
local orthonormal frame field  $\{e_1, \dots , e_n \}$ on $M$,
$\dim M =n$,  $\|f_{*}\|^{2}$ can be expressed as $\|f_{*}\|^{2} =
\langle f_{*}e_{i},f_{*}e_{i}\rangle_N$, where the summation
convention is used. Such a convention will be followed in the
sequel. When a risk of confusion appear, the sum will be written
in detail.

 The critical points of the functional ${\mathcal E}$ on
$C^{\infty}(M,N)$ are known as {\em harmonic maps} and, when $M$
is  closed, they have been characterised as maps with vanishing
{\em tension field}. The tension of $f$ is the vector field
$\tau(f)$ along $f$ which can be locally expressed as
\begin{equation} \label{tensfie}
 \tau(f) = \displaystyle  \widetilde{\nabla}_{e_{i}}(f_{*}e_{i})
  -  f_{*}\nabla_{e_{i}}e_{i},
\end{equation}
where $\widetilde{\nabla}$ denotes the {\em induced connection},
compatible with $\langle \cdot, \cdot \rangle_N$, on the {\em
induced vector bundle} $f^{*}TN = \{ (x,v) \, | \, x \in M, \;
 v \in T_{f(x)}N\}$ of $TN$ by $f$. Here, the fibre on
$x$ in $f^{*}TN$ is identified with $T_{f(x)}N$ and the Equation
\eqref{tensfie} is understood in this way. Then $\tau(f)$ is a
smooth section of $f^{*}TN.$ Denote by $\Gamma^{\infty}(f^{*}TN)$
the space of all smooth sections of $f^{*}TN$, also known as the
space of all {\em variational vector fields} along $f.$ Such a
space can be regarded as the tangent space $T_{f}C^{\infty}(M,N)$
at $f$ of the manifold $C^{\infty}(M,N)$.

If $\bar{N}\subset N$ is a regular submanifold of $N$ such that
$f(M)\subset \bar{N},$ then $f$ belongs to $C^{\infty}(M,\bar{N})$
and the tangent projection ${\rm tan}(V)$ of $V\in
\Gamma^{\infty}(f^{*}TN)$ to $\bar{N}$ is a vector field along
$f:M\to \bar{N}.$ Moreover, for $X\in {\mathfrak X}(M)$ and
$\bar{V}\in \Gamma^{\infty}(f^{*}T\bar{N})$,  we get ${\rm
tan}(\widetilde{\nabla}^{N}_{X}\bar{V}) =
\widetilde{\nabla}_{X}^{\bar{N}}\bar{V},$ where
$\widetilde{\nabla}^{N}$ and $\widetilde{\nabla}^{\bar{N}}$
respectively denote the induced connections via $f:M\to N$ and via
$f:M\to \bar{N}$. Hence, ${\rm tan}\;\tau(f)$ is the tension field
of $f:M\to \bar{N}$ and we have

\begin{lemma}\label{l11} The map
$f:(M, \langle \cdot , \cdot \rangle_M)\to (\bar{N}, \langle \cdot , \cdot \rangle_N)$
is harmonic if and only
if the tension field of $f:(M, \langle \cdot , \cdot \rangle_M)\to
(N, \langle \cdot , \cdot \rangle_N)$ is orthogonal to $\bar{N}.$
\end{lemma}

\section {Harmonicity of sections of sphere  bundles}{\indent} \setcounter{equation}{0}
 Let $\pi:{\mathbb
E}\to M$ be a vector bundle over an $n$-dimensional Riemannian
manifold $(M,\langle \cdot, \cdot \rangle)$ with typical fibre
${\mathbb E}_{x}\cong {\mathbb F}^{m}$, where ${\mathbb F}
=\mathbb R$ or ${\mathbb C}$ and denote by $L({\mathbb E})$ the
principal frame bundle of ${\mathbb E}.$ A point $p$ of
$L({\mathbb E})$ is a pair $(x;p_{1},\dots ,p_{m})$ where $x\in M$
and $\{p_{1},\dots, p_{m}\}$ is a basis of ${\mathbb E}_{x}.$ We
consider a {\em metric connection} $\Gamma: p\in L({\mathbb
E})\mapsto \mathcal{H}_{p}\subset T_{p}L(\mathbb E),$ with respect
to a {\em fibre metric} on ${\mathbb E}$, which we will also
denote by $\langle \cdot, \cdot \rangle$.

For each $v\in {\mathbb E}_{x},$ the corresponding horizontal
subspace $\mathcal{H}_{v}$ of $T_{v}{\mathbb E}$ is defined as
$\mathcal{H}_{v} = (\phi_{\lambda})_{*p}\mathcal{H}_{p},$ where
$(p,\lambda)\in L({\mathbb E})\times {\mathbb F}^{m}$ with $v =
\sum_{i}\lambda^{i}p_{i}$ and $\phi_{\lambda}$ is the mapping
$\phi_{\lambda}: L({\mathbb E}) \to {\mathbb E}$ given by $\phi(q)
= \sum_{i}\lambda^{i}q_{i},$ for all $q\in L({\mathbb E}).$ Hence,
we obtain that $T_{v}{\mathbb E}$ is decomposed into
$T_{v}{\mathbb E} = \mathcal{H}_{v}\oplus \mathcal{V}_{v},$ being
$\mathcal{V}_{v}$ the vertical subspace $\mathcal {V}_{v} =
T_{v}{\mathbb E}_{x},$ and $(M,\langle \cdot,\cdot\rangle)$
acquires a covariant derivative $\nabla$ on the space of the
smooth sections $\Gamma^{\infty}({\mathbb E})$ on ${\mathbb E}$ by
using of the notion of parallel displacement of fibres of
${\mathbb E}.$ Because $\Gamma$ is a metric connection, it follows
that $\nabla$ is {\em compatible with} $\langle \cdot, \cdot
\rangle$, that is, it satisfies
\begin{equation}\label{mc} X \langle \sigma_{1},\sigma_{2} \rangle =
\langle \nabla_{X}\sigma_{1},\sigma_{2} \rangle + \langle
\sigma_{1},\nabla_{X}\sigma_{2}\rangle,
\end{equation}
for all vector field $X$ on $M$ and $\sigma_{1}$, $\sigma_{2}\in
\Gamma^{\infty}({\mathbb E}).$ Moreover, we define $\langle \nabla
\sigma_{1} , \nabla \sigma_{2} \rangle$ in terms of local
orthonormal tangent frames $\{ e_1, \dots, e_n \}$ by the
expression $\langle \nabla \sigma_{1} , \nabla \sigma_{2} \rangle
= \langle \nabla_{e_i} \sigma_{1} , \nabla_{e_i} \sigma_{2}
\rangle$.

The manifold ${\mathbb E}$ admits a Riemannian metric which
generalises the {\em Sasaki metric} of the tangent bundle, see
\cite{Bl}. It will be denoted by $\langle \cdot, \cdot \rangle^S$
and, for $\xi_1,\xi_2\in T_{v}{\mathbb E}$, it is defined by
\begin{equation}\label{metric}
\langle \xi_1,\xi_2\rangle^S =  \langle \pi_* (\xi_1),\pi_*
(\xi_2) \rangle + \langle K(\xi_1),K(\xi_2)\rangle,
\end{equation}
where $K:T{\mathbb E}\to {\mathbb E}$ is the {\em connection map}
\cite{Sak}  of the connection on ${\mathbb E}$ associated with
$\Gamma$. We recall that such a map is  given by $K(\xi) = \iota (
\xi^{\sf v})$, where $\xi^{\sf v}$ is the vertical component of
$\xi\in T{\mathbb E}$ and $\iota$ is the projection $\iota : T
\mathbb E \to \mathbb E$ defined by $\iota(\eta) = 0$ for all
$\eta\in {\mathcal H}$ and $\iota(u_{v}) = u,$ for all $u_{v}\in
\mathcal{V}_{v},$ being $u_{v} = \alpha'(0)$ and $\alpha(t) = v +
tu.$ In particular, if $\sigma$ is a smooth section $\sigma\in
\Gamma^{\infty}({\mathbb E}),$ one obtains that $K(\sigma_{*}X) =
\nabla_{X}\sigma,$ for any vector field $X\in {\mathfrak X}(M).$

If $M$ is compact and oriented, the energy functional of a smooth
section $\sigma\in \Gamma^{\infty}({\mathbb E})$ is defined as the
energy of the map $\sigma:(M,\langle \cdot , \cdot\rangle )\to
({\mathbb E}, \langle \cdot, \cdot \rangle^S).$ Then, from
(\ref{f1}) and using (\ref{metric}), the energy ${\mathcal
E}(\sigma)$ of $\sigma$ can be expressed as
\[
{\mathcal E}(\sigma) = \frac{\textstyle n}{\textstyle 2} {\rm
Vol}(M) + \frac{\textstyle 1}{\textstyle 2}\int_{M}\|\nabla
\sigma\|^{2}dv.
\]
The relevant part of this formula, $B(\sigma) =
\frac{1}{2}\int_{M}\|\nabla \sigma\|^{2}dv$, is usually called the
{\em total bending} of the section $\sigma.$ It is immediate that
$B:\Gamma^{\infty}({\mathbb E})\to \mathbb R$ is always
non-negative and $B(\sigma)$ is zero if and only if $\sigma$ is
$\nabla$-parallel. Thus, the notion of total bending provides a
measure of how a section of $\pi:{\mathbb E}\to M$ fails to be
parallel.

Each smooth section $\sigma\in \Gamma^{\infty}({\mathbb E})$
determines a vertical vector field $\sigma\ve$ on ${\mathbb E}$
given by $\sigma\ve_{v} = \sigma(x)_{v}\in \mathcal{V}_{v},$ for
all $v\in {\mathbb E}_{x},$ and likewise, each vector field $X$ on
$M$ can be lifted to a horizontal vector field $X\h$ on ${\mathbb
E},$ its {\em horizontal lift}.

 In the sequel, we will make use of the {\it musical
isomorphisms} $\flat : TM \to T^* M$ and $\sharp : T^* M \to T M$,
induced by the metric $\langle \cdot , \cdot \rangle$,
respectively defined by $X^{\flat} = \langle X , \cdot \rangle$
and $\langle \theta^{\sharp} , \cdot \rangle = \theta $.

The tension field $\tau(\sigma)$ of $\sigma$ has been
characterised in \cite{GGV} as follows
$$\tau(\sigma) =  \left( (R_{(\sigma,\langle \cdot , \cdot
\rangle)})^{\sharp} \right)\h\comp\sigma - (\nabla^* \nabla \sigma
)\ve\comp \sigma,
$$
where $R_{(\sigma,\langle \cdot , \cdot \rangle)}$ is the one-form
on $M$ given by
\begin{equation} \label{sigerrenab}
R_{(\sigma,\langle \cdot , \cdot \rangle)} (X) = \langle  R_{X, \,
e_i}\sigma, \nabla_{e_i} \sigma \rangle,
\end{equation}
for any vector field $X$ on $M$, where $\{ e_1, \dots , e_n \}$ is
an orthonormal  frame field  and $\nabla^* \nabla \sigma$ is the
{\em connection Laplacian} (or {\em rough Laplacian})\cite{LM}
defined by
$$
\nabla^* \nabla \sigma = -  \left( \nabla^2 \sigma
\right)_{e_i,e_i}.
$$
Here $R_{X,Y}\sigma = \nabla_{[X,Y]}\sigma -
\nabla_X\nabla_Y\sigma + \nabla_Y\nabla_X\sigma$ and
$(\nabla^2\sigma )_{X,Y} = \nabla_X\nabla_Y\sigma -
\nabla_{(\nabla_XY)}\sigma$. Hence, the map $\sigma:(M,\langle
\cdot , \cdot \rangle )\to ({\mathbb E},\langle \cdot ,
\cdot\rangle^S)$ is harmonic if and only if $R_{(\sigma,\langle
\cdot , \cdot \rangle)} = 0$ and $\nabla^* \nabla \sigma = 0$.
Because for $M$ compact, the connection Laplacian of $\sigma$
vanishes if and only if $\sigma$ is parallel (see \cite[page
154]{LM}), it follows that $\sigma$ is harmonic if and only if
$\sigma$ is  parallel.

A critical point $\sigma\in\Gamma^\infty({\mathbb E})$ of the
restriction ${\mathcal E}:\Gamma^{\infty}({\mathbb E}) \to \RR$ of
the energy functional to the space of sections is called a {\em
harmonic section}. Consider $\sigma_{t}\in
\Gamma^{\infty}({\mathbb E})$ a smooth variation of $\sigma$
through sections. Then the corresponding {\em variation vector
field} $x\in M\mapsto V(x) = \frac{d}{dt}_{\mid t=0}\sigma_{t}(x)$
is a section of the induced bundle $\sigma^{*}\mathcal{V}$ of the
vertical subbundle $\sigma^{*}\mathcal{V}\subset T{\mathbb E}$.
Using the first variation formula
\[
\frac{d}{dt}_{\mid t=0}{\mathcal E}(\sigma_{t}) = -\int_{M}\langle
V,\tau(\sigma)\rangle^S dv,
\]
critical points $\sigma\in \Gamma^{\infty}({\mathbb E})$ of the
restriction of ${\mathcal E}$ to $\Gamma^{\infty}({\mathbb E})$
are characterised by the vanishing of the vertical component of
their tension. Hence, we can conclude

\begin{proposition}\label{3o} Let $\pi:{\mathbb E}\to M$ be a vector bundle with
a metric connection over a closed and oriented Riemannian manifold
and $\sigma\in\Gamma^\infty({\mathbb E}).$ Then the following
statements are equivalent:
\begin{enumerate}
\item[{\rm (i)}] the map $\sigma:(M, \langle \cdot , \cdot \rangle
)\to ({\mathbb E}, \langle \cdot , \cdot \rangle^S)$ is harmonic;
\item[{\rm (ii)}] $\sigma$ is a harmonic section; \item[{\rm
(iii)}] $\sigma$ is parallel.
\end{enumerate}
\end{proposition}

Denote by $S_{\mathbb E}(r)$ the sphere  bundle of radius $r>0$ in
${\mathbb E}$ consisting of those elements $v\in {\mathbb E}$ with
$\|v\| = r.$ It is a subbundle and also a hypersurface of
${\mathbb E}.$ For each $\sigma\in \Gamma^{\infty}(S_{\mathbb
E}(r)),$ $\frac{1}{r}\sigma\ve\comp \sigma$ is a unit normal
vector field to $S_{\mathbb E}(r)$ along $\sigma$ and the tangent
projection ${\rm tan}\;\tau(\sigma)$ of $\tau(\sigma)$ to
$S_{\mathbb E}(r)$ is given by
$$
\begin{array}{lcl}
{\rm tan}\;\tau(\sigma) & = & \tau(\sigma) -
\frac{\textstyle 1}{\textstyle r^{2}} \langle \tau(\sigma),\sigma\ve\comp \sigma\rangle^S \sigma\ve\comp \sigma\\[0.5pc]
& = & \left((R_{(\sigma,\langle \cdot , \cdot
\rangle)})^{\sharp}\right)\h\comp\sigma
   - (\nabla^* \nabla \sigma)\ve\comp \sigma
   +   \frac{\textstyle 1}{\textstyle r^{2}}\langle (\nabla^* \nabla \sigma)\ve\comp \sigma,\sigma\ve\comp
   \sigma\rangle^{S}
\sigma\ve\comp \sigma \\[0.5pc]
& = & \left((R_{(\sigma,\langle \cdot , \cdot
\rangle)})^{\sharp}\right)\h\comp \sigma
 + \left( \frac{\textstyle 1}{\textstyle r^{2}} \langle \nabla^* \nabla
\sigma,\sigma \rangle \sigma - \nabla^* \nabla \sigma \right)
 \ve\comp \sigma.
\end{array}
$$
Hence, using Lemma \ref{l11}, we have (see also \cite{GGV})

\begin{proposition}\label{3o1} Let $\pi:{\mathbb E}\to M$ be a vector bundle with
a metric connection over a closed and oriented Riemannian manifold
and $\sigma\in\Gamma^\infty(S_{\mathbb E}(r)).$ Then, we have:
\begin{enumerate}
 \item[{\rm (i)}] the map $\sigma:(M, \langle \cdot , \cdot \rangle
)\to (S_{\mathbb E}(r), \langle \cdot , \cdot \rangle^S)$ is
harmonic if and only if $R_{(\sigma,\langle \cdot , \cdot
\rangle)}  = 0$ and $\nabla^* \nabla \sigma$ is collinear with
$\sigma$.
  \item[{\rm (ii)}] $\sigma$ is a critical point of
${\mathcal E}$ restricted to $\Gamma^\infty(S_{\mathbb E}(r))$ if
and only if $\nabla^* \nabla \sigma$ is collinear with $\sigma$.
\end{enumerate}
\end{proposition}

For general Riemannian manifolds $(M,\langle \cdot ,\cdot
\rangle)$, not necessarily closed and oriented, a section of
$S_{\mathbb E}(r)$ satisfying this last condition is called a {\em
harmonic section of the sphere  bundle}  $S_{\mathbb E}(r)$. If a
harmonic section $\sigma$ is such that $R_{(\sigma,\langle \cdot ,
\cdot \rangle)}=0$, then it is also a harmonic map into
$(S_{\mathbb E}(r), \langle \cdot , \cdot \rangle^S).$ In such a
case, we refer to $\sigma$ as a \emph{harmonic map into a sphere
bundle}.

Let $(\nabla \sigma)^{\mbox{t}} : \Gamma^{\infty}({\mathbb E}) \to
{\mathfrak X}(M)$ be the {\em transpose} operator of $\nabla
\sigma$ with respect to $\langle \cdot , \cdot \rangle$ defined as
\[
\langle(\nabla \sigma)^{\mbox{t}}\varphi,X \rangle = \langle
\varphi,\nabla_{X}\sigma\rangle, \quad \varphi\in
\Gamma^{\infty}({\mathbb E}),\;X\in {\mathfrak X}(M).
\]
The following identity relating  the connection Laplacian and the
{\rm transpose} operator is satisfied \cite[page 155]{LM}
\begin{equation} \label{nanalap}
\langle \nabla^* \nabla \sigma , \varphi \rangle = - \mbox{div} \,
(\nabla \sigma)^{\mbox{t}}\varphi + \langle \nabla \sigma , \nabla
\varphi \rangle,
\end{equation}
which motivates the notation chosen for the connection Laplacian.

If $\|\sigma\|=r$, then (\ref{mc}) implies that $(\nabla
\sigma)^{t}\sigma = 0$. Therefore, using Equation \eqref{nanalap},
we have $ \langle \nabla^* \nabla \sigma,\sigma \rangle = \|
\nabla \sigma \|^{2}$. Hence, $\sigma\in
\Gamma^{\infty}(S_{\mathbb E}(r))$ is a harmonic section of
$S_{\mathbb E}(r)$ if and only if
\begin{equation}
\label{rough}
 \nabla^* \nabla \sigma = \frac{\textstyle 1}{\textstyle
r^{2}}\|\nabla \sigma\|^{2}\sigma.
\end{equation}

The next result will be useful in the discussion of some examples
to decide  whether a harmonic section is a harmonic map.
\begin{lemma} \label{compequi}
Given a harmonic section $\sigma$ of the sphere  bundle
$S_{\mathbb E}(r),$ the one-form $R_{(\sigma,\langle \cdot, \cdot
\rangle)}$ defined in Equation \eqref{sigerrenab} can be also
written as
\begin{equation} \label{generalequis}
R_{(\sigma,\langle \cdot, \cdot \rangle)}(X)  =   \mbox{\rm
div}\left( (\nabla \sigma)^{\mbox{\rm t}}\nabla_{X} \sigma \right)
 +  \langle \nabla_{[X , e_i  ]}
\sigma, \nabla_{e_i} \sigma \rangle - \textstyle \frac12  X \left(
\|\nabla \sigma\|^{2}\right).
\end{equation}
Moreover, if $\langle \nabla_X \sigma , \nabla_Y \sigma \rangle$
is locally expressed by
$$
\langle \nabla_X \sigma , \nabla_Y \sigma \rangle = \sum_{i=1}^n
k_i e_i^{\flat} \otimes e_i^{\flat} (X , Y),
$$
 where
$\{e_1, \dots, e_n\}$ is a local orthonormal frame field and $k_1,
\dots, k_n$ are smooth functions, then
\begin{equation} \label{generalequis1}
R_{(\sigma,\langle \cdot, \cdot \rangle)} = \sum_{i=1}^n
\{e_{i}(k_{i}) + \sum_{j=1}^{n}(k_i-k_j) \langle \nabla_{e_j} e_i,
e_j \rangle \} e^{\flat}_i - \textstyle \frac{1}{2} \displaystyle
d(\sum_{j=1}^{n}k_{j}).
\end{equation}
In particular, if $k_{1} = \dots = k_{n}=\lambda,$ where $\lambda$
is a (non-negative) constant, then $\sigma$ is a harmonic map into
$(S_{\mathbb E}(r),\langle \cdot, \cdot \rangle).$
\end{lemma}
\begin{proof}  From \eqref{sigerrenab}, using the definition of the curvature operator,
 we have
$$
R_{(\sigma,\langle \cdot, \cdot \rangle)}(X) = \langle (\nabla^2
\sigma)_{e_i , X} , \nabla_{e_i} \sigma \rangle - \langle
(\nabla^2 \sigma)_{X, \,e_i \;} , \nabla_{e_i} \sigma \rangle.
$$
$\,$ From this identity, it is immediate to derive
$$
R_{(\sigma,\langle \cdot, \cdot \rangle)}(X) = \langle
\nabla_{e_i} (\nabla_X \sigma) , \nabla_{e_i} \sigma \rangle +
\langle \nabla_{[X,e_i ]} \sigma , \nabla_{e_i} \sigma \rangle-
\langle \nabla_X (\nabla_{e_i} \sigma) , \nabla_{e_i} \sigma
\rangle.
$$
Now,  $\langle \nabla_X (\nabla_{e_i} \sigma) , \nabla_{e_i}
\sigma \rangle =  \frac12 X \left( \|\nabla \sigma\|^{2}\right)$
and, by Equation \eqref{nanalap}, we have
$$
R_{(\sigma,\langle \cdot, \cdot \rangle)}(X)  =   \mbox{\rm
div}\left( (\nabla \sigma)^{\mbox{\rm t}}\nabla_{X} \sigma \right)
 + \langle \nabla^* \nabla
\sigma , \nabla_X \sigma \rangle +
   \langle \nabla_{[X , e_i  ]}
\sigma, \nabla_{e_i} \sigma \rangle - \textstyle \frac12  X \left(
\|\nabla \sigma\|^{2}\right).
$$
Therefore, using \eqref{rough},  Equation \eqref{generalequis}
follows. To show Equation \eqref{generalequis1}, we directly apply
\eqref{generalequis}, taking into account that $ (\nabla
\sigma)^{\mbox{t}} \nabla_{e_i} \sigma = k_i e_i.$
\end{proof}

Finally, we give the first and the second variation formula or the
{\em Hessian form} of the energy functional ${\mathcal E}$
restricted to the set of all sections of the sphere  bundle.
Firstly, we note that $\Gamma^{\infty}({\mathbb E})$ is a module
over the ring of ${\mathbb F}$-valued functions and, for each
$\sigma \in \Gamma^{\infty}(S_{\mathbb E}(r)),$ one obtains the
decomposition  $\Gamma^{\infty}({\mathbb E}) = {\mathcal
V(\sigma)}^{\perp}\oplus {\mathcal V(\sigma)},$ where ${\mathcal
V(\sigma)}$ is the submodule spanned by $\sigma$ and ${\mathcal
V(\sigma)}^{\perp}$ is the orthogonal complement to $\sigma$ on
$\Gamma^{\infty}({\mathbb E}),$ with respect to the metric fibre
in ${\mathbb E}.$ Then $\Gamma^{\infty}(S_{\mathbb E}(r))$ can be
endowed with a structure of Fr\'{e}chet manifold compatible with
its $C^{\infty}$-topology such that each $\sigma\in
\Gamma^{\infty}(S_{\mathbb E}(r))$ is contained in a chart
modelled on ${\mathcal V(\sigma)}^{\perp}$ and, consequently
$T_{\sigma}\Gamma^{\infty}(S_{\mathbb E}(r)) = {\mathcal
V(\sigma)}^{\perp}$. Moreover, a smooth variation $\sigma_{t},$
$t\in ]-\varepsilon, \varepsilon [,$ of $\sigma$ through sections
of $S_{\mathbb E}(r)$ can be regarded as a smooth curve $\gamma:
t\mapsto \gamma(t)= \sigma_{t}$ in $\Gamma^{\infty}(S_{\mathbb
E}(r))$ with $\gamma(0) = \sigma$ and $\gamma'(0)\in {\mathcal
V(\sigma)}^{\perp}$.

\begin{proposition} Let $\pi:{\mathbb E}\to M$ be a vector bundle
with a metric connection over a closed and oriented Riemannian
manifold and let ${\mathcal E}:\Gamma^{\infty}(S_{\mathbb
E}(r))\to \mathbb R$ be the energy functional on
$\Gamma^{\infty}(S_{\mathbb E}(k)).$ We have
\begin{enumerate}
\item [{\rm (i)}] $d{\mathcal E}_{\sigma}(\varphi) =
\displaystyle\int_{M}\langle \nabla^* \nabla \sigma
,\varphi\rangle \, dv,$ for each $\sigma\in
\Gamma^{\infty}(S_{\mathbb E}(r))$ and $\varphi\in {\mathcal
V(\sigma)}^{\perp}$.
 \item [{\rm (ii)}] If $\sigma$ is a  harmonic
section of $S_{\mathbb E}(r)$, then the Hessian form $({\rm
Hess}\;{\mathcal E})_{\sigma}$ on ${\mathcal V(\sigma)}^{\perp}
\cong T_{\sigma}\Gamma^{\infty}(S_{\mathbb E}(r))$ is given by
\[
({\rm Hess}\; {\mathcal E})_{\sigma}\varphi =
\int_{M}(\|\nabla\varphi\|^{2} - \|\varphi\|^{2}\|\nabla
\sigma\|^{2})dv.
\]
\end{enumerate}
\end{proposition}
\begin{proof}
If $\gamma:]-\varepsilon, \varepsilon[ \to
\Gamma^{\infty}(S_{\mathbb E}(r))$ is  a curve such that
$\gamma(0) = \sigma,$ $\gamma'(0) = \varphi\in {\mathcal
V(\sigma)}^{\perp}$, then we obtain

\[
d{\mathcal E}_{\sigma}(\varphi) = \frac{\textstyle d}{\textstyle
dt}_{\mid t=0}{\mathcal E}(\gamma(t)) =
\frac{1}{2}\int_{M}\frac{d}{dt}_{\mid t=0} \| \nabla \gamma \|^2
dv = \int_{M} \langle \nabla \sigma,  \nabla \varphi \rangle dv.
\]
Now, using Equation \eqref{nanalap} and taking into account that
$M$ is closed, we get (i).

For (ii),
$$
\begin{array}{lcl}
({\rm Hess}\;\mathcal E)_{\sigma}\varphi  & = & \frac{\textstyle
d^{2}}{\textstyle dt^{2}}_{\mid t=0}{\mathcal E}(\gamma(t)) =
 \displaystyle \frac12   \int_{M} \frac{\textstyle d^{2}}{\textstyle
dt^{2}}_{\mid t=0}\|\nabla
\gamma\|^{2}dv\\[3mm]
& = & \displaystyle\int_{M}\frac{\textstyle d}{\textstyle
dt}_{\mid t =0} \langle \nabla \gamma , \nabla \gamma'\rangle =
\int_{M} (\|\nabla\varphi\|^{2} + \langle \nabla \sigma
,\nabla\gamma''(0)\rangle)dv.
\end{array}
$$
But, using Equation \eqref{nanalap} as before,  we get
\[
\langle \nabla \sigma, \nabla\gamma''(0)\rangle = {\rm
div}((\nabla\sigma)^{\rm t}\gamma''(0)) + \langle
\gamma''(0),\nabla^* \nabla \sigma \rangle,
\]
and therefore
\[
({\rm Hess}\;\mathcal E)_{\sigma}\varphi =
\int_{M}\left(\|\nabla\varphi\|^{2} + \langle \gamma''(0),\nabla^*
\nabla \sigma \rangle\right) dv.
\]
Now, taking into account   that $\|\gamma(t)\|^{2} = r^2$  and
$\sigma$ is a harmonic section of the sphere bundle, we obtain
\[
\langle \gamma''(0),\nabla^* \nabla \sigma \rangle = \langle
\gamma''(0), \sigma\rangle\|\nabla \sigma\|^{2} = -
\|\varphi\|^{2} \|\nabla\sigma\|^{2}.
\]
Hence (ii) follows.
\end{proof}

\section{Differential forms as  harmonic  maps}
Denote by $\bigwedge^{p}M$ the vector bundle of $p$-forms on $M$
and by $\Omega^{p}M$ the space of its sections, that is, the space
of differential $p$-forms on $M.$ On $\bigwedge^{p}M$ we will
consider the natural fibre metric $\langle \cdot, \cdot \rangle$
defined by
\begin{equation}\label{extendedmetric}
\langle \Psi,\Phi \rangle = \Psi(e_{i_{1}},\dots ,e_{i_{p}}) \Phi
(e_{i_{1}},\dots ,e_{ i_{p}}),
\end{equation}
where $\{e_{1},\dots ,e_{n}\}$ is a local orthonormal frame.
Clearly, the covariant derivative $\nabla$ on $\Omega^{p}M$
obtained as an extension of the Levi Civita connection associated
to the metric $\langle \cdot, \cdot \rangle$ on $M$ is compatible
with such fibre metric, i.e.,  Equation \eqref{mc} is satisfied.

Next Theorem is a first application of Lemma \ref{compequi} and
will be extremely useful in working with examples to be able to
claim that certain harmonic  sections of  some particular sphere
 bundles are also harmonic maps.
\begin{theorem} \label{thpair}
Let $(M,\langle \cdot, \cdot \rangle)$ be an $n$-dimensional
Riemannian manifold and $(\Psi,\Phi)$ a pair of differential forms
of constant length $\|\Psi\|= r_{1}$ and $\|\Phi\| = r_{2},$ $\Psi
\in \Omega^p M $ and $\Phi \in \Omega^{p+1} M .$ If $\nabla_X \Psi
= \lambda X \lrcorner \Phi$ and $\nabla_X \Phi = \mu X^{\flat}
\wedge \Psi$, where $\lambda,$ $\mu$ are constants and $0 \leq p <
n$, then $\Psi$ and $\Phi$ are harmonic maps into the
corresponding sphere bundles $S_{\Omega^p M }(r_{1})$ and
$S_{\Omega^{p+1} M }(r_{2}).$
\end{theorem}
\begin{proof} For $x \in M$, taking an orthonormal frame field $\{
e_1, \dots, e_n \}$, such that $\left( \nabla_{e_i} e_j \right)_x
=0$, it follows that
 \begin{gather*}
 \left(
\nabla^{\ast} \nabla \Psi \right)_x = - \nabla_{e_{ix}}
\left(\nabla_{e_i} \Psi \right) = -  \lambda\mu e_{ix} \lrcorner (
e_{ix}^{\flat} \wedge \Psi_x) = -(n-p) \lambda\mu \Psi_x , \\
\left( \nabla^{\ast} \nabla \Phi \right)_x = -  \nabla_{e_{ix}}
\left(\nabla_{e_i} \Phi \right) =  - \lambda\mu e_{ix}^{\flat}
\wedge ( e_{ix} \lrcorner \Phi_x) = -(p+1)\lambda\mu \Phi_x.
 \end{gather*}
 Hence, $(\Psi,\Phi)$ is a pair of harmonic sections of the respective sphere  bundles and, using \eqref{rough},  we
 get
 \[
 \|\nabla \Psi\|^{2} = -(n-p)\lambda\mu r_{1}^{2},\;\;\;\;\; \|\nabla
 \Phi\|^{2} = -(p+1)\lambda\mu r_{2}^{2}.
 \]
 Next, we compute $R_{(\Psi ,\langle \cdot, \cdot \rangle)}$ and
 $R_{(\Phi ,\langle \cdot, \cdot \rangle)}.$
 It is straightforwardly
obtained
$$
 \mbox{div} \left((\nabla \Psi)^{\mbox{t}} \nabla_X \Psi\right)) = \textstyle \frac12 (n-p+1)
\lambda\mu X(\| \Psi \|^2) +  \lambda^2 \langle e_i \lrcorner
\Phi, \nabla_{e_i} X \lrcorner \Phi \rangle = \lambda^2 \langle
e_i \lrcorner \Phi, \nabla_{e_i} X \lrcorner \Phi \rangle.
$$
Also it is  direct to obtain
\begin{equation*}
\begin{split}
\langle \nabla_{[X,\cdot]} \Psi , \nabla_{\cdot} \Psi \rangle = &
\lambda^2\{  \langle \nabla_X e_i \lrcorner \Phi , e_i \lrcorner
\Phi \rangle -  \langle e_i \lrcorner \Phi, \nabla_{e_i} X
\lrcorner
\Phi \rangle \}\\
= & - \lambda^{2}\{\textstyle\frac{1}{2(p+1)} X (\|\Phi\|^2) +
\langle e_i \lrcorner \Phi, \nabla_{e_i} X \lrcorner \Phi
\rangle\}
\\
=& - \lambda^2 \langle e_i \lrcorner \Phi, \nabla_{e_i} X
\lrcorner \Phi \rangle.
\end{split}
\end{equation*}
 Therefore, by Lemma \ref{compequi}, we have
 $R_{(\Psi ,\langle \cdot, \cdot \rangle)}=0$. We also obtain
 \begin{equation*}
\begin{split}
 \mbox{div} \left(  (\nabla \Phi)^{\mbox{t}}\nabla_X \Phi\right) =&  \textstyle \frac12 (p+1)
\mu X(\lambda\|\Phi \|^2 + \mu \| \Psi \|^2) + \mu^{2} \langle
e_i^{\flat} \wedge  \Psi,
(\nabla_{e_i} X)^{\flat} \wedge \Psi \rangle\\
 =& \mu^{2} \langle e_i^{\flat}
\wedge \Psi, (\nabla_{e_i} X)^{\flat} \wedge \Psi \rangle.
\end{split}
\end{equation*}
On the other hand, we have
\begin{equation*}
\begin{split}
\langle \nabla_{[X,e_i]} \Phi , \nabla_{e_i} \Phi \rangle = &
\mu^{2}\{\langle (\nabla_X e_i)^{\flat} \wedge  \Psi , e_i^{\flat}
\wedge \Psi \rangle - \langle ( \nabla_{e_i} X)^{\flat} \wedge
\Psi , e_i^{\flat} \wedge \Psi \rangle \}\\
= & \mu^{2}\{\textstyle \frac12 (p+1)(2p+1-n) X (\|\Psi\|^2)  -
\langle e_i^{\flat} \wedge \Psi, (\nabla_{e_i} X)^{\flat} \wedge
\Psi \rangle\}
\\
=& - \mu^2 \langle e_i^{\flat} \wedge \Psi, (\nabla_{e_i}
X)^{\flat} \wedge \Psi \rangle.
\end{split}
\end{equation*}
Now, using again Lemma \ref{compequi}, we have
 $R_{(\Phi ,\langle \cdot, \cdot \rangle)}=0$ and then $(\Psi,\Phi)$
is moreover a pair of harmonic  maps into their respective sphere
  bundles.
\end{proof}

\section{Examples of harmonic maps}

First we recall  some notions relative to $\Lie{G}$-structures,
where $\Lie{G}$ is a subgroup of the linear group $\Lie{GL}(n ,
\mathbb R)$. An $n$-dimensional manifold $M$ is equipped with a
$\Lie{G}$-structure, if its frame bundle admits a reduction to the
subgroup $\Lie{G}$. If $M$ possesses a $\Lie{G}$-structure, then
there always exists a $\Lie{G}$-connection defined on $M$.
Moreover, if $(M ,\langle \cdot , \cdot \rangle)$ is an oriented
Riemannian $n$-manifold and  $\Lie{G}$ is a closed and connected
subgroup of $\SO(n)$, then there exists a unique metric
$\Lie{G}$-connection $\nabla^G = \nabla + \xi^{\Lie{G}}$ such that
$\xi^{\Lie{G}} \in  T^* M \otimes \lie{g}^{\perp}$, where
$\lie{g}^{\perp}$ denotes the orthogonal complement in $\so(n)$ of
the Lie algebra $\lie{g}$ of $\Lie{G}$ and $\nabla$ denotes the
Levi-Civita connection \cite{CleytonSwann:torsion}. The tensor
$\xi^{\Lie{G}}$ is called the {\it intrinsic torsion}  of the
$\Lie{G}$-structure and $\nabla^{\Lie{G}}$ is said to be the {\it
minimal $\Lie{G}$-connection}.

\subsection{Nearly K{\"a}hler 6-manifolds} An almost Hermitian
manifold is a  Riemannian $2n$-manifold $(M, \langle \cdot , \cdot
\rangle)$ endowed with an almost complex structure $J$ compatible
with the metric. The presence of such a structure is equivalent to
say that $M$ is equipped with a $\Lie{U}(n)$-structure. Under the
action of $\Lie{U}(n)$, the space $T^* M \otimes
\lie{u}(n)^{\perp}$ of possible intrinsic torsion tensors
$\xi^{\Lie{U}(n)}$  is decomposed into irreducible
$\Lie{U}(n)$-modules:
\begin{enumerate}
 \item if $n=1$, $ \xi^{\Lie{U}(1)} \in T^* M \otimes \un(1)^\perp
= \{ 0 \}$;
  \item if $n=2$, $ \xi^{\Lie{U}(2)} \in T^* M \otimes
\un(2)^\perp = \Wc2^{\Lie{U}(2)} +
  \Wc4^{\Lie{U}(2)}$;
\item if $n \geqslant 3$, $ \xi^{\Lie{U}(n)} \in T^* M \otimes
\un(n)^\perp =
  \Wc1^{\Lie{U}(n)} + \Wc2^{\Lie{U}(n)} + \Wc3^{\Lie{U}(n)} + \Wc4^{\Lie{U}(n)}$.
\end{enumerate}
where $\mathcal W_i^{\Lie{U}(n)}$ are the irreducible
$\Lie{U}(n)$-modules  given  by Gray and Hervella
\cite{Gray-H:16}. The vanishing of $\Un(n)$-components of
$\xi^{\Lie{U}(n)}$ gives rise to a natural classification of
almost Hermitian manifolds. Associated with the almost Hermitian
structure, it is usually considered the two form $\omega = \langle
\cdot , J \cdot \rangle$, called the {\it K{\"a}hler form}. One
can use the $\Lie{U}(n)$-isomorphism $\xi^{\Lie{U}(n)} \to -
\xi^{\Lie{U}(n)} \omega = \nabla \omega$ to identify the intrinsic
$\Lie{U}(n)$-torsion with $\nabla \omega$. Thus,  Gray and
Hervella showed   conditions expressed by means of $\nabla \omega$
to characterise classes of almost Hermitian manifolds.

 Each fibre $T_x M$ of the
tangent bundle  can be consider as a complex vector space by
defining $i v = Jv$.  We will write $T_x M_{\mathbb C}$ when we
are regarding $T_x M$ as such a space. If, for all $x \in M$,
there exists a complex volume $n$-form on  $T_x M_{\mathbb C}$
defined  by
$$
\Psi_x = \left( \Psi_{+} \right)_{x} + i
\left(\Psi_{-}\right)_{x},
$$
such that $\Psi_+$ and $\Psi_-$ are real global differential
$n$-forms on $M$ compatible with the almost Hermitian structure,
then  $M$ is said to be a {\it special almost Hermitian manifold}
(see \cite{Cabrera:special} for details). Such a fact is
equivalent to say that $M$ is equipped with an
$\Lie{SU}(n)$-structure. For higher dimensions, $n \geq 4$, the
space $T^* M \otimes \lie{su}(n)^{\perp}$ is decomposed into five
irreducible $\Lie{SU}(n)$-modules $\mathcal W_1^{\Lie{SU}(n)},
\dots , \mathcal W_5^{\Lie{SU}(n)}$. The first four modules are
such that $\mathcal W_i^{\Lie{SU}(n)}= \mathcal W_i^{\Lie{U}(n)}$,
$i=1, \dots, 4,$ and  $\mathcal W_5^{\Lie{SU}(n)} \cong T^* M$.

For $n=3$, the space $T^* M \otimes \lie{su}(3)^{\perp}$ of
intrinsic $\Lie{SU}(3)$-torsion tensors is decomposed into  the
following modules (\cite{CS,Cabrera:special})
$$
T^* M \otimes \lie{su}(3)^{\perp} = \mathcal W_1^+ + \mathcal
W_1^- + \mathcal W_2^+ + \mathcal W_2^- + \mathcal
W_3^{\Lie{SU}(3)} + \mathcal W_4^{\Lie{SU}(3)} + \mathcal
W_5^{\Lie{SU}(3)},
$$
where $\mathcal W_i^+ + \mathcal W_i^- = \mathcal
W_i^{\Lie{U}(3)}$, $i=1,2$,   $\mathcal W_j^{\Lie{SU}(3)} =
\mathcal W_j^{\Lie{U}(3)}$, $j=3,4$,  $\mathcal W_5^{\Lie{SU}(3)}
\cong T^* M$,  $\mathcal W_1^+ \cong \mathcal W_1^- \cong  \mathbb
R$ and $\mathcal W_2^+ \cong \mathcal W_2^- \cong  \lie{su}(3)$.

 When
$\xi^{\Lie{U}(n)} \in \mathcal W_1^{\Lie{U}(n)}$, the almost
Hermitian manifold  is called nearly K{\"a}hler.
 Gray \cite{Gray:nearly-Kaehler} showed that any nearly K{\"a}hler and non-K{\"a}hler  6-manifold is Einstein.
 Furthermore, the Einstein constant
$\rho$ is positive. In this  case, one can consider the 3-form
$\Psi_+$ of type $(3,0)$ such that $3w_1^+ \Psi_+ = d \omega$,
where $5 \left( w_1^+ \right)^2 = \rho$. Now we define
$$
\Psi_- = - \Psi_+ ( J \cdot , \cdot , \cdot)
$$
and  fix $ \Psi_+ + i \Psi_-$ as complex volume form, obtaining an
$\Lie{SU}(3)$-structure (special almost Hermitian structure) on
the manifold. Such an $\Lie{SU}(3)$-structure is of type $\mathcal
W_1^+$ (in this case the $\mathcal W^{\Lie{SU}(3)}_5$-part of the
intrinsic torsion vanishes, see \cite{Cabrera:special}). Thus, we
will have a two-form $\omega$ and a three-form $\Psi_+$, such that
$$
\nabla_X \omega = w_1^+  X \lrcorner \Psi_+, \qquad  \nabla_{X}
\Psi_+ = - w_1^+ X^{\flat} \wedge \omega, \qquad \| \omega \|^2  =
6, \qquad \| \Psi_+ \|^2 = 24.
$$
Hence, the conditions contained  in Theorem \ref{thpair} are
satisfied. Then we get
$$
\nabla^* \nabla \omega =  4 \left( w_1^+ \right)^2 \omega, \quad
\nabla^* \nabla \Psi_+ =  3 \left( w_1^+ \right)^2 \Psi_+, \quad
\quad R_{(\omega, \langle \cdot ,\cdot \rangle)} = R_{( \Psi_+ ,
\langle \cdot ,\cdot \rangle)}  =0.
$$

Additionally, we can also  consider the pair consisting of the
three-form $\Psi_-$ and the four-form $\omega \wedge \omega$. Such
forms satisfy
\begin{gather*}
\nabla_{X} \Psi_- =  \nabla_{JX} \Psi_+ = -  w_1^+ JX^{\flat}
\wedge \omega = \tfrac12 w_1^+ X \lrcorner (\omega \wedge
\omega),\\
\nabla_{X} \left( \omega \wedge \omega \right) =   2  w_1^+ \left(
X \lrcorner \Psi_+ \right) \wedge \omega =  - 2  w_1^+ X^{\flat}
\wedge  \Psi_-, \\
\| \Psi_- \|^2 = 24, \qquad \| \omega \wedge \omega \|^2  = 144.
\end{gather*}
Now, making use again of Theorem \ref{thpair}, we have
$$
\nabla^* \nabla \Psi_- =  3 \left( w_1^+ \right)^2 \Psi_-, \quad
\nabla^* \nabla (\omega \wedge \omega) =  4 \left( w_1^+ \right)^2
\omega \wedge \omega, \quad \quad R_{(\Psi_-,\langle \cdot , \cdot
\rangle)} = R_{(\omega \wedge \omega,\langle \cdot , \cdot
\rangle)} =0.
$$
In summary,
\begin{theorem} For a nearly K\"{a}hler $6$-manifold, the
differential forms $\omega$, $\Psi_+$, $ \Psi_-$ and $ \omega
\wedge \omega$ are  harmonic  maps into their respective sphere
bundles.
\end{theorem}

\subsection{Nearly parallel $\Lie{G}_2$-manifolds}  A Riemannian
seven-manifold $M$ admits a $\Lie{G}_2$-structu\-re if and only if
there exists a  three-form $\phi$ on $M$,  nowhere zero, which is
$\Lie{G}_2$-invariant and it is locally expressed by
$$
\phi = \textstyle  \sum_{i \in \mathbb Z_7} e_i^{\flat} \wedge
e_{i+1}^{\flat} \wedge e_{i+3}^{\flat} ,
$$
where $\{ e_0 , \dots, e_6 \}$ are certain local orthonormal frame
fields. Such frames, usually called {\it Cayley frames}, are
adapted to the $\Lie{G}_2$-structure and the seven-form
$e_0^{\flat} \wedge \dots \wedge e_6^{\flat}  = Vol$ is globally
defined and fixed as volume form. Thus, corresponding to the
volume form there is a Hodge star operator $\ast$. The four-form
$\ast \phi$ is also $\Lie{G}_2$-invariant and locally expressed by
$$
\ast \phi = - \textstyle \sum_{i \in \mathbb Z_7} e_{i+2}^{\flat}
\wedge e_{i+4}^{\flat}  \wedge e_{i+5}^{\flat}  \wedge
e_{i+6}^{\flat}.
$$

Associated to the $\Lie{G}_2$-structure, we have the minimal
$\Lie{G}_2$-connection $\nabla^{\Lie{G}_2} = \nabla +
\xi^{\Lie{G}_2}$, such that $\xi^{\Lie{G}_2} \in T^* M \otimes
\lie{g}_2^{\perp} \subset T^*M \otimes \Lambda^2 T^* M$. In this
case, the action of $\Lie{G}_2$ decomposes  the space $T^* M
\otimes \lie{g}_2^{\perp}$ of possible intrinsic torsion tensors
into four $\Lie{G}_2$-irreducible modules $\mathcal
W_1^{\Lie{G}_2}, \dots , \mathcal W_4^{\Lie{G}_2}$
\cite{FernandezGray}. If one considers the $\Lie{G}_2$-modules of
bilinear forms on $\mathbb R^7$ equipped with the usual Euclidean
product $\langle \cdot, \cdot \rangle$, one can show  that
$\mathcal W_1^{\Lie{G}_2} \cong \mathbb R \langle \cdot, \cdot
\rangle$, $ \mathcal W_2^{\Lie{G}_2} \cong \lie{g}_2 $,   $
\mathcal W_3^{\Lie{G}_2} \cong S_0^2 \mathbb R^{7*}$, $ \mathcal
W_4^{\Lie{G}_2} \cong \lie{g}_2^{\perp} \cong \mathbb R^7 $.
 When $\xi^{\Lie{G}_2} \in \mathcal
 W_1^{\Lie{G}_2}$, the $\Lie{G}_2$-structure is called {\it nearly
 parallel}.
 In such a case, the manifold is Einstein \cite{Gr3},
  $ \nabla_X \phi = \tfrac{k}{4} X \lrcorner \ast \phi$, $ \nabla_X \ast \phi = - \tfrac{k}{4} X^{\flat}
  \wedge  \phi$ and $ \rho = \frac{k^2}{16} $  is the Einstein constant
  \cite{Cabrera:g2}. Since  $ 4 \| \phi \|^2 = \| \ast \phi \|^2 = 7 . 4! $, we are in the conditions
  of Theorem \ref{thpair}. Therefore,
$$
\nabla^* \nabla \phi =    \tfrac{k^2}{4} \phi =  4 \rho \phi,
\quad \nabla^* \nabla \ast \phi =   \tfrac{k^2}{4} \ast \phi =  4
\rho \ast \phi, \quad R_{(\phi,\langle \cdot, \cdot \rangle)} =
R_{(\ast \phi,\langle \cdot, \cdot \rangle)} =0.
$$
In summary,
\begin{theorem} For a nearly parallel $\Lie{G}_2$-manifold, the
differential forms $\phi$ and $\ast \phi$ are harmonic maps into
their respective sphere   bundles.
\end{theorem}

\subsection{a-Sasakian manifolds} An almost contact metric manifold
is a Riemannian $(2n+1)$-manifold $(M,\langle \cdot , \cdot
\rangle)$ equipped with a $(1,1)$-tensor $\varphi$ and a one-form
$\eta$ such that $\langle \varphi X , \varphi Y \rangle = \langle
X , Y \rangle - \eta(X) \eta(Y)$ and $\varphi^2 = - I + \eta
\otimes \zeta$, where $\zeta^{\flat} = \eta$. The presence of the
mentioned tensors on the manifold is equivalent to say that $M$ is
endowed  with a $\Lie{U}(n) \times 1$-structure, where $\Lie{U}(n)
\times 1$ is considered such that $\Lie{U}(n) \times 1 \subseteq
\Lie{SO}(2n+1)$. In this case, the cotangent space at each point $
T^*_x M$ is not irreducible under the action of the group
$\Lie{U}(n) \times 1$. In fact, $T^*_x M = \mathbb R \eta +
\eta^{\perp}$ and
$$
\lie{so}(2n+1) \cong \Lambda^2 T^* M \cong \Lambda^2 \eta^{\perp}
+ \eta^{\perp} \wedge \mathbb R \eta = \lie{u}(n) +
\lie{u}(n)^{\perp}_{| \zeta^{\perp}} + \eta^{\perp} \wedge \mathbb
R \eta.
$$
Therefore, for the space $T^* M \otimes \lie{u}(n)^{\perp}$ of
possible intrinsic $\Lie{U}(n) \times 1$-torsion, we  have
$$
T^* M \otimes  \lie{u}(n)^{\perp} = \eta^{\perp} \otimes
\lie{u}(n)^{\perp}_{| \zeta^{\perp}} \, + \,  \eta \otimes
\lie{u}(n)^{\perp}_{| \zeta^{\perp}} \, + \, \eta^{\perp} \otimes
\eta^{\perp} \wedge \eta  \, + \, \eta \otimes
  \eta^{\perp} \wedge   \eta.
$$
Chinea  and Gonz{\'a}lez-D{\'a}vila \cite{ChineaGonzalezDavila}
showed that $T^* M \otimes  \lie{u}(n)^{\perp}$ is decomposed into
twelve irreducible $\Lie{U}(n)$-modules $\mathcal C_1, \dots ,
\mathcal C_{12}$, where
\begin{eqnarray*}
\eta^{\perp} \otimes \lie{u}(n)^{\perp}_{| \zeta^{\perp}}  & = &
\mathcal
C_1 + \mathcal C_2 + \mathcal C_3 + \mathcal C_4, \\
\eta^{\perp} \otimes \eta^{\perp} \wedge   \eta  & = & \mathcal
C_5 + \mathcal C_8 + \mathcal C_9 + \mathcal C_6 +
\mathcal C_7  + \mathcal C_{10}, \\
 \eta \otimes \lie{u}(n)^{\perp}_{| \zeta^{\perp}} & = & \mathcal C_{11} , \\
 \eta \otimes \eta^{\perp} \wedge  \eta & = & \mathcal
C_{12}.
\end{eqnarray*}
Note that $\mathcal C_1, \dots , \mathcal C_4$ are the Gray and
Hervella's $\Lie{U}(n)$-modules, i.e.,  $\mathcal{C}_i \cong
\mathcal W_i^{\Lie{U}(n)}$. Furthermore, note that $\varphi$
restricted to $\zeta^{\perp}$ works as an almost complex structure
and, if one considers the $\Lie{U}(n)$-action on the bilinear
forms $\bigotimes^2 \eta^{\perp}$, then we have the decomposition
$$
\textstyle \bigotimes^2 \eta^{\perp} = \mathbb R \langle \cdot ,
\cdot \rangle_{|\zeta^{\perp}} + \lie{su}(n)_s  + \sigma^{2,0} +
\mathbb R F  + \lie{su}(n)_a  +
\lie{u}(n)^{\perp}_{\zeta^{|\perp}},
$$
where $F$ is the form, called the {\it fundamental two-form},
defined by $F = \langle \cdot , \varphi \cdot \rangle$.

The modules  $\lie{su}(n)_s$ (resp., $\lie{su}(n)_a$) consist of
Hermitian symmetric (resp., skew-symmetric) bilinear forms
orthogonal to $\langle \cdot , \cdot \rangle_{|\zeta^{\perp}}$
(resp., $F$)
 and  $\real{\sigma^{2,0}}$ ($\lie{u}(n)^{\perp}_{|\zeta^{\perp}}
 $) is the space of  anti-Hermitian
symmetric (resp., skew-symmetric) bilinear forms. With respect to
the modules $\mathcal C_i$, one has $\eta^{\perp} \otimes
\eta^{\perp} \wedge  \mathbb R \eta \cong   \bigotimes^2
\eta^{\perp}$ and, using the $\Lie{U}(n)$-map $\xi^{\Lie{U}(n)}
\to - \xi^{\Lie{U}(n)} \eta = \nabla \eta$, it is obtained
$$
\mathcal C_5 \cong  \mathbb R \langle \cdot , \cdot
\rangle_{|\zeta^{\perp}}, \quad \mathcal C_8 \cong \lie{su}(n)_s,
\quad \mathcal C_9 \cong \real{\sigma^{2,0}}, \quad \mathcal C_6
\cong \mathbb R F, \quad \mathcal C_7 \cong  \lie{su}(n)_a, \quad
\mathcal C_{10} \cong \lie{u}(n)^{\perp}_{|\zeta^{\perp}}.
$$

In summary,   under the action of $U(n)\times 1$, the space of
possible intrinsic torsion tensors
  $\mbox{T}^* M\otimes  \lie{u}(n)^{\perp}$ is decomposed
into:
\begin{enumerate}
\item if $n=1$, $ \xi^{\Lie{U}(1)} \in \mbox{T}^* M \otimes
\un(1)^\perp = \mathcal C_{5} + \mathcal C_{6} + \mathcal C_{9} +
\mathcal C_{12}$; \item if $n=2$, $ \xi^{\Lie{U}(2)} \in
\mbox{T}^* M \otimes \un(2)^\perp = \mathcal C_{2} + \mathcal
C_{4} + \dots + \mathcal C_{12}$; \item if $n \geqslant 3$, $
\xi^{\Lie{U}(n)} \in \mbox{T}^* M \otimes \un(n)^\perp =
  \mathcal C_{1}  + \dots  + \mathcal C_{12}$.
\end{enumerate}

Here, we will consider the particular case that the intrinsic
$\Lie{U}(n) \times 1$-torsion is contained in $\mathcal C_6$. In
such a case, $M$ is called an $a$-Sasakian manifold (Sasakian when
a=1) and it is characterised by the condition
$$
\nabla_X F = - a  X^{\flat} \wedge \eta.
$$
Note that $\mbox{\it d}^* F (\zeta) = 2n a$, where $d^*$ means the
coderivative. Moreover, if $n \geq 2$, then $a$ is constant
\cite{Marrero}.

Since for almost contact metric manifolds one has $\left(\nabla_X
\eta \right) Y = \left( \nabla_X F \right) (\zeta , \varphi Y)$,
then, for an $a$-Sasakian manifold, it follows $ \nabla_X \eta = a
X \lrcorner F$.

Now, we consider $a$-Sasakian $(2n+1)$-manifolds, $n\geq 2$, and
the pairs of differential forms $(\eta \wedge F^r , F^{r+1} )$,
where $0 \leq  r\leq n$ and $F^r = F \wedge \stackrel{(r)}{\dots}
\wedge F$. Such pairs satisfy
$$
\nabla_X \left( \eta \wedge F^r \right) =  \tfrac{1}{r+1} \, a \,
X  \lrcorner F^{r+1}, \qquad \nabla_X  F^{r+1} = -(r+1)a \,
X^{\flat} \wedge \eta \wedge  F^r.
$$
Thus, they satisfy the conditions required in Theorem
\ref{thpair}. Note that $\| \eta \wedge F^r\|^2 = \frac{(2r+1)! r!
n!}{(n-r)!}$ and $\| F^{r+1}\|^2 = \frac{(2(r+1))!
(r+1)!n!}{(n-r-1)!}$. Therefore, we get
\begin{gather*}
\nabla^* \nabla \left( \eta \wedge F^r \right) =  2(n-r)a^2 \eta
\wedge F^r , \quad \nabla^* \nabla F^{r+1} =  2(r+1)a^2 F^{r+1},
\\
R_{(\eta \wedge F^r,\langle \cdot, \cdot \rangle)} = R_{(
F^{r+1},\langle \cdot, \cdot \rangle)}  =  0.
\end{gather*}

In summary,
\begin{theorem} For an $a$-Sasakian $(2n+1)$-manifold, the
differential forms $\eta \wedge F^r$ and $F^{r+1}$, $0 \leq r \leq
n$,  are   harmonic  maps into their respective sphere bundles.
\end{theorem}
Note that, in general, $a$-Sasakian manifolds are not Einstein.
Diverse interesting  aspects of the  so-called Sasakian-Einstein
manifolds can be found  in
\cite{BoyerGalicki:Sasakian-Einstein,BoyerGalicki:5Sasakian-Einstein}.

\subsection{3-a-Sasakian manifolds} A $(4n+3)$-manifolds $M$ is
said to be endowed with an {\it almost contact metric
3-structure}, if $M$ has a Riemannian metric $\langle \cdot ,
\cdot \rangle$ and three almost contact metric structures
$(\varphi_i, \zeta_i , \eta_i)$, $i=1,2,3$, satisfying
\begin{gather*}
\eta_i ( \zeta_j) = \delta_{ij}, \quad \varphi_i (\zeta_j) = -
\varphi_j (\zeta_i) = \zeta_k, \quad \eta_i \circ \varphi_j = -
\eta_j \circ \varphi_i = \eta_k, \\
\varphi_i \circ \varphi_j - \eta_j \otimes \zeta_i = - \varphi_j
\circ \varphi_i - \eta_i \otimes \zeta_j = \varphi_k, \qquad
\langle \varphi_i X , \varphi_i Y \rangle = \langle X , Y \rangle
- \eta_i(X) \eta_i (Y),
\end{gather*}
for all cyclic permutation $(ijk)$ of $(123)$. In the language of
$\Lie{G}$-structures, one says  that the manifold is equipped with
an $\Lie{Sp}(n) \times I_3$-structure. Here, it is considered
$\Lie{Sp}(n) \times I_3 \subseteq \Lie{SO}(4n) \times I_3
\subseteq \Lie{SO}(4n+3)$.

Now, we will consider $(4n+3)$-manifolds endowed  with an almost
contact metric 3-structure such that the three structures are
$a$-Sasakian. If $a_i$ is the constant corresponding to the
structure $i$, $i=1,2,3$, it is immediate to show that
$a_1=a_2=a_3=a$. In this case,  the manifold is called a {\it
$3$-$a$-Sasakian manifold} ({\it $3$-Sasakian}, when $a=1$).
Kashiwada  \cite{Kas} proved that 3-Sasakian manifolds are
Einstein, being $\rho = 2(2n+1)$ the Einstein constant.

On $3$-$a$-Sasakian manifolds,  we consider the pair of
differential forms $(\Psi^{(r)}, \Omega^{(r)})$, where
$$
 \Psi^{(r)} =
\textstyle \sum_{i=1}^3  \eta_i \wedge  F_i^r, \qquad \Omega^{(r)}
= \textstyle \sum_{i=1}^3 F_i^{r+1}.
$$
The tensors $\Omega$ and $\Psi$ are of constant length and
$$
\nabla_X \Psi^{(r)} = -  \textstyle \frac{a}{r+1} X \lrcorner
\Omega^{(r)}, \qquad \nabla_X \Omega^{(r)} = (r+1) a X \wedge
\Psi^{(r)}.
$$
Now, making use of Theorem \ref{thpair}, we have
$$
\nabla^* \nabla \Psi^{(r)} =  2(2n+1-r) a^2 \Psi^{(r)}, \quad
\nabla^* \nabla \Omega^{(r)} = 2 (r+1) a^2 \Omega^{(r)}, \quad
R_{( \Psi^{(r)}, \langle \cdot , \cdot \rangle)} =
R_{(\Omega^{(r)}, \langle \cdot,\cdot \rangle)} =0.
$$

Another pair of differential forms to be considered consists  of
$\textstyle \cyclic_{ijk}^{123} \eta_i \wedge F_j \wedge F_k \; $
and $\; F_1 \wedge F_2 \wedge  F_3$, where  $\cyclic$ denotes
cyclic sum over the listed elements.
 For such forms, we have
$$
\textstyle \nabla_X \left( \cyclic_{ijk}^{123} \eta_i \wedge F_j
\wedge F_k \right) = - a X \lrcorner \left( F_1 \wedge F_2 \wedge
F_3 \right)
$$
and
$$
\textstyle \nabla_X \left(  F_1 \wedge F_2 \wedge F_3 \right) = a
X \wedge \left(\cyclic_{ijk}^{123} \eta_i \wedge F_j \wedge F_k
\right).
$$
Since these forms are of constant length, we can use Theorem
\ref{thpair}. Therefore,
\begin{gather*}
\textstyle \nabla^* \nabla  \left( \cyclic_{ijk}^{123} \eta_i
\wedge F_j \wedge F_k \right) =  2(2n-1) a^2 \left(
\cyclic_{ijk}^{123} \eta_i \wedge F_j \wedge F_k \right)  ,\\
\textstyle \nabla^* \nabla  \left( F_1 \wedge F_2 \wedge F_3
\right) =  6 a^2 F_1 \wedge F_2 \wedge F_3 , \\
 \textstyle R_{( \cyclic_{ijk}^{123} \eta_i \wedge F_j \wedge F_k
 , \langle \cdot , \cdot \rangle)} = R_{( F_1 \wedge F_2 \wedge
 F_3, \langle \cdot , \cdot \rangle)} =0.
\end{gather*}

Now, we consider the pair $\eta_i \wedge F_j + \eta_j \wedge F_i$
and $F_i \wedge F_j$, where $i\neq j$. Since these forms are of
constant length and
$$
\nabla_X \left( \eta_i \wedge F_j + \eta_j \wedge F_i \right) = -
a X \lrcorner (F_i \wedge F_j), \qquad \nabla_X \left( F_i \wedge
F_j \right) = a X \wedge \left( \eta_i \wedge F_j + \eta_j \wedge
F_i \right),
$$
then, by Theorem \ref{thpair},
$$
\nabla^* \nabla \left( \eta_i \wedge F_j + \eta_j \wedge F_i
\right) =  4n a^2 \left( \eta_i \wedge F_j + \eta_j \wedge F_i
\right), \qquad \nabla^* \nabla (F_i \wedge F_j) =  4 a^2 F_i
\wedge F_j,
$$
$$
R_{(\eta_i \wedge F_j + \eta_j \wedge F_i ,\langle \cdot, \cdot
\rangle)} = R_{( F_i \wedge F_j,\langle \cdot, \cdot \rangle)} =
0.
$$

In summary,
\begin{theorem} For a 3-$a$-Sasakian $(4n+3)$-manifold, the
differential forms $\sum_{i=1}^3  \eta_i \wedge F_i^r$,
$\sum_{i=1}^3 F_i^{r+1}$, $\cyclic_{ijk}^{123} \eta_i \wedge F_j
\wedge F_k$, $F_1 \wedge F_2 \wedge  F_3$, $\eta_i \wedge F_j +
\eta_j \wedge F_i$ and $F_i \wedge F_j$, where  $0 \leq r \leq
2n+1$ and $(ijk)$ is a cyclic permutation of $(123)$, are harmonic
 maps into their respective sphere  bundles.
\end{theorem}

Finally, for 3-$a$-Sasakian manifolds, we will also find some
differential forms which are eigenvectors with respect to the
connection Laplacian and they do not follow the scheme contained
in Theorem \ref{thpair}. In fact, we will discuss the harmonicity
as a map of  such forms  by applying Lemma \ref{compequi}. For
instance, let us consider the three-form
$$
\textstyle \vartheta = (2n+3) \eta_1 \wedge \eta_2 \wedge  \eta_3
+ \sum_{i=1}^3 \eta_i \wedge F_i.
$$
If $\{e_1, \dots ,e_{4n},e_{4n+1}, e_{4n+2}, e_{4n+3}\}$ is a
local orthonormal frame field such that  at $x \in M$ is the basis
for vectors $\{ e_{1 \, x}, \dots , e_{4n \, x}, e_{4n+1 \, x} =
\zeta_{1 \, x},
 e_{4n+2 \, x} =  \zeta_{2 \, x} ,  e_{4n+3 \, x} =  \zeta_{3 \, x} \}$
 adapted to the 3-Sasakian structure
and $\left( \nabla_{e_i} e_j \right)_x=0$, then
$$
\textstyle \nabla^* \nabla \eta_1 \wedge \eta_2 \wedge \eta_3 = 12
(n+1) a^2 \eta_1 \wedge \eta_2 \wedge  \eta_3 + 4 a^2 \sum_{i=1}^3
\eta_i \wedge F_i  .
$$
Now, from this last identity, $\nabla^* \nabla \Psi^{(1)} = 4n
\Psi^{(1)}$ and using the identities
\begin{equation} \label{equrem}
\textstyle  \sum_{r=1}^{4n+3} (e_r \lrcorner (e_r \wedge \eta_i)=
2(2n+1) \eta_i, \qquad \sum_{r=1}^{4n+3} (e_r \lrcorner F_i)
\wedge (e_r \lrcorner F_j) = -2 F_k - \eta_i \wedge \eta_j,
\end{equation}
a straightforward computation shows  that $\nabla^* \nabla
\vartheta  =  12  (n+1) a^2 \vartheta$.
 \vspace{2mm}

Now we  compute $R_{(\vartheta, \langle \cdot, \cdot \rangle)}$.
One can deduce that $\vartheta$ is of constant length and
$$
\langle \nabla_{X} \vartheta , \nabla_Y \vartheta \rangle = 18
(4n^2 +10 n+3) a^2 \langle X , Y \rangle - 6 (12n^2+26n+9) a^2
\textstyle \sum_{i=1}^3 \eta_i(X) \eta_i(Y).
$$
Now, making use of Lemma \ref{compequi}, we have
$$
\textstyle R_{(\vartheta, \langle \cdot, \cdot \rangle)} = 6
(12n^2+26n +9) a^2 \sum_{i=1}^{4n} \sum_{j=1}^3 \left(  \eta_j(
\nabla_{\zeta_j} e_i) e^{\flat}_i - \langle \nabla_{e_i} \zeta_j ,
e_i \rangle \eta_j \right).
$$
But, for $3$-$a$-Sasakian manifolds, one obtains
$$
\textstyle \nabla \zeta_i = a \sum_{s=1}^{n} \left(   \varphi_i
e_s^{\flat}  \otimes e_s - e_s^{\flat} \otimes  \varphi_i e_s +
\varphi_j e_s^{\flat}  \otimes \varphi_k e_s  - \varphi_k
e_s^{\flat}  \otimes \varphi_j e_s  \right) +   a ( \eta_k \otimes
\zeta_j -  \eta_j \otimes \zeta_k),
$$
where $(ijk)$ is a cyclic permutation of $(123)$. Therefore, for
$i=1, \dots , 4n$ and $j=1,2,3$, it follows
\begin{equation} \label{escalprod}
\eta_j( \nabla_{\zeta_j} e_i)  =  \langle \nabla_{e_i} \zeta_j ,
e_i \rangle =0,
\end{equation}
and, as a consequence, $R_{(\vartheta, \langle \cdot , \cdot
\rangle)} = 0$.
 \vspace{2mm}

In the same context of $3$-$a$-Sasakian manifolds, we have
 \begin{eqnarray*}
 \textstyle \nabla_X \left( \cyclic_{ijk}^{123} \eta_i \wedge \eta_j \wedge
F_k \right)
 &= & \textstyle a \cyclic_{ijk}^{123}   \eta_i \wedge  (X \lrcorner F_j) \wedge
 F_k   - a \cyclic_{ijk}^{123}  \eta_i \wedge F_j \wedge  (X \lrcorner F_k) \\
 &&   \textstyle  -  3 a X^{\flat} \wedge \eta_1 \wedge \eta_2 \wedge
 \eta_3.
 \end{eqnarray*}
From this identity we get
 \begin{eqnarray*}
 \textstyle  \nabla_{e_{s \, p}} \left( \nabla_{e_s}
  \left( \cyclic_{ijk}^{123} \eta_i \wedge \eta_j \wedge
F_k \right) \right)
 &= & \textstyle  2 a^2 \cyclic_{ijk}^{123}   (e_s \lrcorner F_i) \wedge  (e_s \lrcorner F_j) \wedge
 F_k
 \\[-2mm]
 && \textstyle - 2 a^2 \cyclic_{ijk}^{123}  \eta_i  \wedge  e_s \lrcorner (e_s \wedge \eta_j) \wedge
 F_k   \\
 &&  \textstyle  - 5 a^2 \cyclic_{ijk}^{123}  \eta_i  \wedge \eta_j \wedge e_s \wedge  (e_s \lrcorner F_k).
 \end{eqnarray*}
 Thus, taking Equations  \eqref{equrem} into account,
 we obtain
\begin{eqnarray*}
  \textstyle  \nabla^*  \nabla  \left( \cyclic_{ijk}^{123} \eta_i \wedge \eta_j \wedge
F_k \right)  &=& \textstyle 8 (n+2) a^2 \cyclic_{ijk}^{123} \eta_i
\wedge \eta_j \wedge F_k  + 4 a^2 \Omega^{(1)}.
\end{eqnarray*}
Now, since $\nabla^* \nabla \Omega^{(1)} = 4 \Omega^{(1)}$, it
follows
$$
 \textstyle  \nabla^*  \nabla  \left( \Omega^{(1)} + (2n+3) \cyclic_{ijk}^{123} \eta_i \wedge \eta_j \wedge
F_k \right) = 8(n+2)a^2 \left( \Omega^{(1)} + (2n+3)
\cyclic_{ijk}^{123} \eta_i \wedge \eta_j \wedge F_k \right).
$$
Writing $\Psi =  \Omega^{(1)} + (2n+3) \cyclic_{ijk}^{123} \eta_i
\wedge \eta_j \wedge F_k $, it is straightforward to check
$$
\langle \nabla_X \Psi , \nabla_Y \Psi \rangle = 6 (4n+1) a^2
\langle X , Y \rangle - 6 a^2 \textstyle \sum_{i=1}^3 \eta_i(X)
\eta_i(Y).
$$
Now, making use of Lemma \ref{compequi}, we have
$$
R_{(\Psi, \langle \cdot , \cdot \rangle)} = \textstyle
 6 a^2 \sum_{i=1}^{4n} \sum_{j=1}^3 \left(
\eta_j( \nabla_{\zeta_j} e_i) e^{\flat}_i - \langle \nabla_{e_i}
\zeta_j , e_i \rangle \eta_j \right).
$$
Again, using Equation \eqref{escalprod} and the fact that $\Psi$
is of constant length, we obtain $R_{(\Psi, \langle \cdot , \cdot
\rangle)} =0$.
 \vspace{2mm}

Finally, since for 3-$a$-Sasakian manifolds  we have
\begin{equation*}
\begin{split}
\nabla_X (\eta_i \wedge \eta_j) = \,& a (X\lrcorner F_i) \wedge
\eta_j
+ \eta_i \wedge (X \lrcorner F_j), \\
 \nabla_{e_{s\, p}} \left( \nabla_{e_s} (\eta_i \wedge \eta_j)
\right)  =\,& a^2 ( 2 (e_s \lrcorner F_i) \wedge (e_s \lrcorner
F_j) - 2 \eta_i \wedge \eta_j + \eta_i(e_s) e_s^{\flat} \wedge
\eta_j + \eta_j(e_s) \eta_i \wedge e_s^{\flat} ),
\end{split}
\end{equation*}
we obtain $ \nabla^* \nabla (\eta_i \wedge \eta_j) = 2 a^2( 2 F_k
+ (4n+3) \eta_i \wedge \eta_j)$. As a consequence, for any $(ijk)$
cyclic permutation, we have
$$
\nabla^* \nabla ( F_k + (2n+1)\eta_i \wedge \eta_j) = 2 (4n+3)
a^2(  F_k + (2n+1) \eta_i \wedge \eta_j).
$$
Moreover, it is direct to check that
\begin{equation*}
\begin{split}
 \langle \nabla_X ( F_k +
(2n+1)\eta_i \wedge \eta_j) , \nabla_Y ( F_k + (2n+1)\eta_i \wedge
\eta_j) \rangle  = &2(2(2n+1)^2 +1) a^2 \langle X,Y \rangle \\
  & \hspace{-4mm} -
2(2n^2+2n+1)a^2 \textstyle \sum_{i=1}^3 \eta_i(X)\eta_i(Y).
\end{split}
\end{equation*}
Therefore, using Lemma \ref{compequi} and Equation
\eqref{escalprod}, we get
$$
R_{(F_k + (2n+1)\eta_i \wedge \eta_j, \langle \cdot , \cdot
\rangle)} = \textstyle 2(2n^2+2n+1) a^2 \sum_{i=1}^{4n}
\sum_{j=1}^3 \left( \eta_j ( \nabla_{\zeta_j} e_i) e^{\flat}_i -
\langle \nabla_{e_i} \zeta_j , e_i \rangle \eta_j \right)=0.
$$
In summary,
\begin{theorem} For a 3-$a$-Sasakian $(4n+3)$-manifold, the
differential forms
\begin{equation*}
\begin{split}
&\textstyle (2n+3) \eta_1 \wedge \eta_2 \wedge  \eta_3 +
\sum_{i=1}^3
\eta_i \wedge F_i,\\
&\textstyle \sum_{i=1}^3 F_i \wedge F_i + (2n+3)
\cyclic_{ijk}^{123} \eta_i
\wedge \eta_j \wedge F_k,\\
&F_k + (2n+1)\eta_i \wedge \eta_j,
 \end{split}
 \end{equation*}
 where  $(ijk)$ is a cyclic
permutation of $(123)$, are   harmonic  maps into their respective
sphere  bundles.
\end{theorem}

\subsection{b-Kenmotsu manifolds} Next we will consider almost
contact metric manifolds $M$ such that the intrinsic $\Lie{U}(n)
\times 1$-torsion is contained in $\mathcal C_5$. In such a case,
$M$ is known  as a \emph{$b$-Kenmotsu manifold} and it is
characterised by the condition
$$
\nabla_X F = b \eta \wedge (X \lrcorner F).
$$
Note that $\nabla_X \eta = - b X^{\flat} + b \eta (X) \eta$ and
$\mbox{div} \, \zeta = - d^* \eta = -2 n b$, where $d^*$ stands
for  the coderivative as we have pointed out above. The exterior
derivative will be denoted by $d$.
 For b-Kenmotsu manifolds, $\eta$ is closed and $b$ is a function such that $d b = f \eta$
 \cite{Marrero}.
\begin{proposition} For $b$-Kenmotsu manifolds, we have:
\begin{eqnarray*}
 \nabla^* \nabla (F^r) & = & 2 r b^2 F^r,  \\
 R_{(F^r, \langle \cdot , \cdot \rangle)} & = & 2r b \| F^r \|^2(rb^2-f) \eta,
 \\
 \nabla^* \nabla (\eta \wedge F^r) & = &2(n-r) b^2 \eta \wedge
 F^r,  \\
 R_{(\eta \wedge F^r, \langle \cdot , \cdot \rangle)} & = & 2 (2r+1)(n-r) \|F^r\|^2 b(b^2-
f) \eta,
\end{eqnarray*}
where $0 \leq r \leq n$. Moreover, we also have:
\begin{enumerate}
\item[{\rm (a)}] $R_{(F^r, \langle \cdot ,\cdot \rangle)}=0$ if
and only if one of the following equivalent conditions is
satisfied:
\begin{enumerate}
 \item[{\rm (i)}] $\mbox{\it d}b(\zeta)=rb^2$,
 \item[{\rm (ii)}]$2n \, \mbox{\rm grad} \left( \mbox{\rm div}\,(\zeta) \right) =
- r \mbox{\rm div}^2\,(\zeta) \zeta$,
 \item[{\rm (iii)}] $2n\, \Delta \eta = r (\mbox{\it
d}^* \eta)^2 \eta$;
 \end{enumerate}
\item[{\rm (b)}]
 $ R_{(\eta \wedge F^r, \langle
\cdot , \cdot \rangle)}=0$ if and only if one of the following
equivalent conditions is satisfied:
\begin{enumerate}
 \item[{\rm (i)}] $\mbox{\it d}b(\zeta)=b^2$,
 \item[{\rm (ii)}]$2n \, \mbox{\rm grad} \left( \mbox{\rm div}\,(\zeta) \right) =
- \mbox{\rm div}^2\,(\zeta) \zeta$,
 \item[{\rm (iii)}]$2n\, \Delta \eta = (\mbox{\it
d}^* \eta)^2 \eta$;
 \end{enumerate}
\end{enumerate}
 where $\Delta$ denotes the Hodge Laplacian, $\Delta = \mbox{\it dd}^* +\mbox{\it d}^* \mbox{\it d}$.
\end{proposition}
Note that $R_{(F^r, \langle \cdot ,\cdot \rangle)}$ and  $R_{(\eta
\wedge F^r, \langle \cdot , \cdot \rangle)}$ are always closed
and, in particular, if $\mbox{\rm div} \, \zeta = - \mbox{\it d}^*
\eta = -2n b$ is a non-zero constant, then $R_{(F^r, \langle \cdot
,\cdot \rangle)}$ and $R_{(\eta \wedge F^r, \langle \cdot , \cdot
\rangle)}$ are nowhere zero.
\begin{proof}
 By straightforward  computation we obtain  the expressions for
 $\nabla^* \nabla (F^r)$ and $\nabla^* \nabla (\eta \wedge F^r)$.
Furthermore, from the identities
$$
\langle \nabla_X F^r , \nabla_Y F^r \rangle =  \textstyle
 \frac{r^2 }{n}  b^2\| F^r \|^2 ( \langle X,Y \rangle - \eta(X)
\eta(Y))
$$
and
$$
\langle \nabla_X (\eta \wedge F^r) , \nabla_Y (\eta \wedge F^r)
\rangle = \textstyle  \frac{(2r+1)(n-r)}{n} b^2 \| F^r \|^2
(\langle X, Y \rangle - \eta (X) \eta(Y) ),
$$
using Lemma \ref{compequi}, the expressions for  $R_{(F^r, \langle
\cdot , \cdot \rangle)}$ and $R_{(\eta \wedge F^r, \langle \cdot ,
\cdot \rangle)}$ follow. From these last expressions the remaining
parts of Proposition are immediately  deduced.
\end{proof}

\begin{example}{\rm
The one-form $\eta$ is closed for $b$-Kenmotsu manifolds. This
implies that one has local coordinates $(x_1, \dots ,\- x_{2n},
x_{2n+1})$ such that $d x_{2n+1} = \eta$. If $k$ is a constant, we
have:
\begin{enumerate}
\item[(a)] $F^r$ is a harmonic map if and only if
$$
b (x_1, \dots, x_{2n} , x_{2n+1}) = -\frac{1}{r(x_{2n+1} +k)},
$$
\item[(b)]  $\eta \wedge F^r$ is a harmonic map if and only if
$$
b (x_1, \dots, x_{2n} , x_{2n+1}) = -\frac{1}{x_{2n+1} +k}.
$$
\end{enumerate}

 In \cite{Marrero}, Marrero considers  the product manifold $M=L\times V$,
  where $L$ is $S^1$ or an open interval, and
 $(V,J,\langle \cdot,\cdot \rangle_V)$ is a K\"{a}hler $2n$-manifold.
 Let $E$ be a nowhere vanishing vector field on $L$ and $\sigma$ a
 positive function on $L$. Taking
 $$
 \varphi(cE,X) = (0,JX), \quad \zeta = (E,0), \quad
 \langle (cE,X), (dE,Y)\rangle = cd + \sigma \langle X, Y \rangle_V,
 $$
where $c$,$d$ are functions on $L$ and  $X,Y$ are vector fields on
$V$,  a $b$-Kenmotsu structure is obtained, $ b=-\frac{1}{2} d
(\mbox{ln} \sigma) (E)$. In particular, we  consider $L$ as an
open interval and $E=\frac{\partial}{\partial t}$, where $t$ is
the coordinate. For $K$ and $C \neq 0$ constants, two interesting
particular cases are:
\begin{enumerate}
\item[(a)$'$]
 If we choose $\sigma = C(t+K)^{\frac{2}{r}}$,   we will obtain that $F^r $ is a
harmonic  map into sphere   bundle. In fact, one  checks  that $b=
- \frac{1}{r(t+K)}$ and $\mbox{\it d}b(\zeta)=r b^2=
\frac{1}{r(t+K)^2}$.

\item[(b)$'$] If we choose $\sigma = C(t+K)^2$,  we will obtain
that $\eta \wedge F^r $ is a harmonic map into a sphere bundle. In
fact, one obtains $b= - \frac{1}{t+K}$ and $\mbox{\it
d}b(\zeta)=b^2= \frac{1}{(t+K)^2}$.
\end{enumerate}
 }
\end{example}

\subsection{Locally conformal parallel $p$-forms}{\indent}
Throughout this subsection, we are assuming that $\Psi$ is a
$p$-form of constant length.
\begin{definition} {\rm
A $p$-form $\Psi$ on a Riemannian $n$-manifold $M$ is said to be
{\it locally conformal parallel}, if there exists a closed
one-form $\theta$ on $M$ such that
\begin{equation} \label{lcpdef}
\nabla_X \Psi = X^\flat \wedge (\theta^\sharp \lrcorner \Psi) -
\theta \wedge (X \lrcorner \Psi),
\end{equation}
for all $X \in \mathfrak X (M)$.  We will refer to the one-form
$\theta$ as the Lee form of $\Psi$. }
\end{definition}
The following results will be useful in examples.
 \begin{proposition} \label{lcpgeneral}
 If $\Psi$ is a locally conformal parallel $p$-form on a Riemannian
$n$-manifold $M$ with Lee form $\theta$, then its coderivative
$d^* \Psi$ and its rough Laplacian $\nabla^* \nabla \Psi$ are
respectively given by
 \begin{gather*}
 d^* \Psi = (p-n) \theta^\sharp \lrcorner \Psi,
\\
 \nabla^* \nabla \Psi = p \| \theta \|^2 \Psi + (n-2p) \theta
 \wedge (\theta^\sharp \lrcorner \Psi).
 \end{gather*}
 In particular, if $2p=n$, then $\Psi$ is a harmonic
section of its  corresponding sphere  bundle.
\end{proposition}
\begin{proof} The expression for $d^* \Psi$ is obtained by a
direct computation. In order to compute $(\nabla^* \nabla
\Psi)_m$, for $m \in M$, we will consider a local orthonormal
frame field $\{ e_1 , \ldots , e_n \}$ such that $(\nabla_{e_i}
e_j )_m =0$. Thus, because $\Psi$ is locally conformal parallel,
we have
$$
(\nabla^* \nabla \Psi)_m = - e_i^\flat \wedge  (\theta^\sharp
\lrcorner (\nabla_{e_i} \Psi)) - e_i^\flat \wedge ( (\nabla_{e_i}
\theta)^\sharp \lrcorner \Psi) + \nabla_{e_i} \theta \wedge (e_i
\lrcorner \Psi) - \theta \wedge d^* \Psi.
$$
Now, using the expression for $\nabla \Psi$ given in
\eqref{lcpdef} and the identities $e_i \wedge (e_i \lrcorner \Psi)
= p \Psi$ and $e_i^\flat \wedge \theta \wedge (\theta^\sharp
\lrcorner ( e_i \lrcorner \Psi))= (p-1) \theta \wedge
(\theta^\sharp \lrcorner \Psi)$, we obtain
$$
 - e_i^\flat \wedge  (\theta^\sharp
\lrcorner (\nabla_{e_i} \Psi)) = - p \theta \wedge ( \theta^\sharp
\lrcorner \Psi) + p \| \theta \|^2 \Psi.
$$
Moreover, because $\theta$ is closed, we have $(\nabla_X
\theta)(Y) = (\nabla_Y \theta)(X)$ and  it is not hard to see that
$$
 e_i^\flat \wedge ( (\nabla_{e_i}
\theta)^\sharp \lrcorner \Psi) = \nabla_{e_i} \theta \wedge (e_i
\lrcorner \Psi).
$$
Finally, from all  this,  the  required expression for $\nabla^*
\nabla \Psi$ follows.
\end{proof}

In the examples below, in order to apply Lemma \ref{compequi}, we
will need to compute $\langle \nabla_X \Psi, \nabla_Y \Psi
\rangle$. For doing this we will use \eqref{lcpdef} to obtaine the
expression
\begin{eqnarray} \label{scnaxy}
\langle \nabla_X \Psi, \nabla_Y \Psi \rangle  & = & p \|
\theta^\sharp \lrcorner \Psi \|^2 \langle X , Y \rangle + p \|
\theta \|^2
\langle X \lrcorner \Psi  , Y \lrcorner \Psi \rangle \\
&&
  - p \theta(X) \langle \theta^\sharp \lrcorner \Psi  , Y \lrcorner \Psi \rangle
  - p \theta(Y) \langle \theta^\sharp \lrcorner \Psi  , X \lrcorner \Psi \rangle
   \nonumber
  \\
  &&
  - 2 p (p-1) \langle X \lrcorner ( \theta^\sharp \lrcorner \Psi)  , Y \lrcorner (\theta^\sharp \lrcorner \Psi)
  \rangle. \nonumber
\end{eqnarray}

\begin{example}[Locally conformal K\"{a}hler manifolds]{\rm An
almost Hermitian $2n$-manifold $(M, \langle \cdot , \cdot \rangle,
J)$ of type $\mathcal W_4$ is characterised  by the condition that
its K\"{a}hler two-form is given by
$$
\nabla_X \omega = X^\flat \wedge (\theta^\sharp \lrcorner \omega)
- \theta \wedge (X \lrcorner \omega).
$$
$\,$From this identity and the nature of $\omega$,  one can deduce
that $\theta$ has to  be closed. Thus, $\omega$ is a locally
conformal parallel two-form.
\begin{proposition} \label{lcpgeneral1}
If the  Lee form $\theta$ of $\omega$ is not zero somewhere and $r
< n$, then $\omega^r$ is a harmonic  section of its corresponding
sphere bundle if and only if $2r = n$.
\end{proposition}
\begin{proof} Direct computations show that
\begin{equation} \label{scpr1}
 \|\omega^r \|^2 = \textstyle  \frac{(2r)! r! n!}{(n-r)!}, \qquad
 \| X \lrcorner \omega^r \|^2 = \textstyle  \frac{(2r)!  r! n!}{2n(n-r)!} \| X
\|^2.
\end{equation}
Moreover, since $\nabla_X \omega^r = X^\flat \wedge (\theta^\sharp
\lrcorner \omega^r) - \theta \wedge (X \lrcorner \omega^r)$, it
follows from Proposition \ref{lcpgeneral} that
$$
 \nabla^* \nabla \omega^r = 2r  \| \theta \|^2 \omega^r + 2(n-2r) \theta
 \wedge (\theta^\sharp \lrcorner \omega^r).
 $$
Additionally, if $\nabla^* \nabla \omega^r = f \omega^r$,  taking
into account that  $\langle X^\flat \wedge \beta , \gamma \rangle
= p \langle \beta , X \lrcorner \gamma \rangle$  into account for
any $(p-1)$-form $\beta$ and  $p$-form $\gamma$,  we have
$$
 f \| \omega^r \|^2 =  2r  \| \theta \|^2 \| \omega^r \|^2 + 4r(n-2r) \| \theta^\sharp \lrcorner
 \omega^r \|^2.
$$
Now, making use of \eqref{scpr1}, we obtain
\begin{equation} \label{firstvalue:f}
  f   =   \textstyle  \frac{4r (n-r)}{n}  \| \theta \|^2.
\end{equation}
On the other hand, if we consider a local adapted  orthonormal
frame $\{ e_1 , J e_1 , \dots , e_n , J e_n \}$ such that $\theta
= \| \theta \| e_1$, by  computing $\nabla^* \nabla \omega^r ( e_1
, Je_1 , \ldots , e_r , J e_r)$, we have
$$
 2r \| \theta \|^2 r! + 2(n-2r) \| \theta \|^2 r!  = f r!.
$$
Therefore, $f= 2 (n-r) \| \theta \|^2$. From this and
\eqref{firstvalue:f}, it follows $n=2r$.
\end{proof}

In order to compute $\langle \nabla_X \omega^r , \nabla_Y \omega^r
\rangle$, a straightforward computation shows that
\begin{eqnarray*}
 \langle X
\lrcorner ( \theta \lrcorner \omega^r ), Y \lrcorner ( \theta
\lrcorner \omega^r) \rangle & = &    \textstyle   \frac{r  r!
(2r-2)! (n-2)!}{(n-r)!}
    \left\{  (r-1) \|\theta\|^2   \langle  X , Y \rangle \right.
   \\
 &&
  \left.   - (r-1)\theta(X) \theta(Y) + (n-r)  J \theta(X) J
  \theta(Y)\right\}.
\end{eqnarray*}
From this,  \eqref{scnaxy} and \eqref{scpr1}, we obtain
\begin{eqnarray*}
 \langle \nabla_X \omega^r , \nabla_Y \omega^r
\rangle &= &
             \textstyle  \frac{2r \,  r!  (2r)!(n-2)!  }{(n-r-1)!} ( \| \theta \|^2     \langle X,Y \rangle
               - \theta (X) \theta (Y)
               - J \theta (X) J \theta (Y)).
\end{eqnarray*}
Finally,  choosing a local orthonormal frame field  as above such
that $\theta = \| \theta \| e_1$, we will have $ \| \theta \|
\mbox{div}(e_1) = - d^* \theta - d(\| \theta \|)(e_1)$ and  $  \|
\theta \| \mbox{div}(Je_1) = - d^* (J\theta) - d(\| \theta
\|)(Je_1)$. Moreover, since $\theta$ is closed, then
$\nabla_{\theta^\sharp} \theta = \| \theta \| d (\| \theta \|)$.
Taking all  this into account  and the fact that $2r=n$ as well,
and making use of Lemma \ref{compequi}, we will obtain
$$
R_{(\omega^r,\langle \cdot , \cdot \rangle)}
 = \textstyle  \frac12 (n-1)
(n!)^2 \left( -  (n-2)  d (\| \theta \|^2)
   + d^* \theta \,\theta +  d^* (J\theta) \,J\theta + \nabla_{J
   \theta^\sharp} J \theta \right).
$$
In particular, if the Lee form $\theta$ is parallel, then it is
not hard to see that $\nabla_{J
   \theta^\sharp} J \theta = 0$ and   $d^* (J \theta) = d^* \omega (\theta^\sharp) =
2 (n-1) J \theta (\theta^\sharp) = 0$. Therefore,
$R_{(\omega^r,\langle \cdot , \cdot \rangle)}=0$.
\begin{theorem} \label{lckparlee}
For a locally conformal K\"{a}hler $4n$-manifold, $\omega^n$ is a
harmonic map into its sphere  bundle if and only if
$$
  (n-2)  d (\| \theta \|^2) =  (d^* \theta) \,\theta +  d^* (J\theta) \,J\theta + \nabla_{J
   \theta^\sharp} J \theta,
$$
where $\theta$ is the Lee form of $\omega$. In particular, if
$\theta$ is parallel, then $\omega^n$ is a harmonic map into its
sphere   bundle.
\end{theorem}
Hopf manifolds are diffeomorphic to $S^1 \times S^{2m-1}$ and
admit a locally conformal K\"{a}hler structure with parallel Lee
form $\theta$ \cite{Va}. Furthermore, $\theta$ is nowhere zero. In
particular, if $m=2n$, then  $S^1 \times S^{4n-1}$ satisfies the
conditions of last Theorem. In general, if we consider a Sasakian
$(4n-1)$-manifold $M$, then the product manifold $ \mathbb R
\times M$ can be equipped with a locally conformal K\"{a}hler
structure satisfying the conditions of Theorem \ref{lckparlee}. }
\end{example}

\begin{example}[Locally conformal parallel
$\Lie{Spin}(7)$-manifolds] {\rm

 Now, let us consider $\mathbb R^8$
endowed with an orientation and its standard inner product. Let
$\{e,e_0,...,e_6\}$ be an oriented orthonormal basis.  Consider
the four-form $\Phi$ on $\mathbb R^8$ given by
\begin{eqnarray}\label{1}
\Phi &=& \sum_{i \in \mathbb{Z}_7}  e\wedge e_i \wedge e_{i+1}
\wedge  e_{i+3} - \sigma \sum_{i \in \mathbb{Z}_7} e_{i+2} \wedge
e_{i+4} \wedge e_{i+5} \wedge  e_{i+6},
\end{eqnarray}
where $\sigma$ is a fixed constant such that $\sigma = + 1$ or
$\sigma=-1$, and $+$ in the subindexes means the sum in
$\mathbb{Z}_7$.  We fix $e \wedge e_0 \wedge \dots \wedge e_6 =
\frac{\sigma}{14} \Phi \wedge \Phi$ as a volume form.

The subgroup of $\GL(8,\mathbb R)$ which fixes $\Phi$ is
isomorphic to the double covering $\Spin(7)$ of $\SO(7)$
\cite{HL}. Moreover, $\Spin(7)$ is a compact simply-connected Lie
group of dimension 21 \cite{Br}. The Lie algebra $\lie{spin}(7)$
of $\Spin(7)$ is isomorphic to the skew-symmetric two-forms $\psi$
satisfying the  linear equations
$$
\sigma \psi ( e_{i},e) + \psi ( e_{i+1}, e_{i+3}) + \psi (
e_{i+4}, e_{i+5}) + \psi ( e_{i+2}, e_{i+6}) = 0,
$$
for all $i \in \mathbb{Z}_7$. Shortly, $\lie{spin}(7)\cong  \{\psi
\in \Lambda^2 T*M |*( \psi\wedge\Phi)= \psi \}$, where $*$ is the
Hodge star operator. The orthogonal complement
$\lie{spin}(7)^{\perp}$ of $\lie{spin}(7)$ in $\Lambda^{2} \mathbb
R^{8\ast} = \lie{so}(8)$ is the seven-dimensional space generated
by
\begin{equation}\label{new1}
\beta_i = \sigma e_{i} \wedge e + e_{i+1} \wedge e_{i+3} + e_{i+4}
\wedge e_{i+5} + e_{i+2} \wedge e_{i+6},
\end{equation}
where $i \in \mathbb{Z}_7$. Equivalently, $\lie{spin}(7)^{\perp}$
is described  as the space  consisting of those skew-symmetric
two-forms $\psi$ such that $*( \psi\wedge\Phi)= - 3 \psi$.

A $\Spin(7)$\emph{-structure} on an eight-manifold $M^8$ is by
definition a reduction of the structure group of the frame bundle
to $\Spin(7)$; we shall also say that $M$ is a
\emph{$\Spin(7)$-manifold}. This can be geometrically described by
saying that there exists a nowhere vanishing global differential
four-form $\Phi$ on $M^8$ and a local frame $\{e,e_0,\dots,e_6\}$
such that the four-form $\Phi$  can be locally written as in
\eqref{1}. The four-form $\Phi$ is called the \emph{fundamental
form} of the $\Spin(7)$-manifold $M$ \cite{Bo} and the  local
  frame  $\{e, e_0, \dots, e_6 \}$ is called {\it a Cayley frame}.

The fundamental form of a $\Spin(7)$-manifold determines a
Riemannian metric $\langle \cdot , \cdot\rangle$ through
 $\langle X,Y \rangle = - \frac17 \ast \left( (X \lrcorner \Phi)
 \wedge \ast \left( Y \lrcorner \Phi \right) \right)$
 \cite{Gr}. Thus, $\langle \cdot , \cdot \rangle$ is called the metric induced by $\Phi$.
 Any  Cayley frame becomes an orthonormal frame with respect to such a metric.

There are four classes of $\Spin(7)$-manifolds given by
Fern{\'a}ndez in \cite{F}. They are obtained as irreducible
$\Spin(7)$-representations   of the space $\overline W \cong
\mathbb R^{8*} \otimes \lie{spin}(7)^{\perp}$ of all possible
covariant derivatives $\nabla\Phi$.   The Lee form $\theta$ is
defined by
\begin{equation}\label{c2}
\theta = -\textstyle  \frac{1}{7}*(*d\Phi\wedge\Phi)= \textstyle
\frac17*(d^* \Phi\wedge \Phi).
\end{equation}
Alternatively, the classification can be described in terms of the
Lee form as follows : $\overline W_0 : d\Phi=0; \quad \overline
W_1 : \theta=0; \quad \overline W_2 : d\Phi = \theta \wedge\Phi;
\quad \overline W :\overline W=\overline W_1+ \overline W_2$
\cite{C1}.

A $\Spin(7)$-structure of the class $\overline W_2$ is called
 a \emph{ locally conformal parallel $\Spin(7)$-structure}.
 This is motivated by the fact that   the Lee form of a
$\Spin(7)$-structure in the class $\overline W_2$ is closed.
Therefore,  such a manifold is locally conformal  to a parallel
$\Spin(7)$-manifold. Furthermore, the covariant derivative $\nabla
\Phi$ for such manifolds is given by
\begin{equation} \label{covderspin}
 4 \nabla_X \Phi =      X^{\flat} \wedge (\theta^{\sharp}
\lrcorner \Phi) -   \theta \wedge (X \lrcorner \Phi).
\end{equation}
\begin{theorem} \label{cfult}
The fundamental form $\Phi$ of a locally conformal parallel
$\Lie{Spin}(7)$-structure is a  harmonic section of its
corresponding sphere bundle. Furthermore, if $\theta$ denotes the
Lee form of the $\Lie{Spin}(7)$-structure, then $\Phi$ is a
harmonic  map into its sphere  bundle if and only if $(d^* \theta)
\, \theta = 3 d \left( \| \theta \|^2 \right)$.
\end{theorem}
In particular, since the $\Lie{Spin}(7)$-structure
 defined  on the product of spheres $S^1 \times S^7$ in \cite{C1}  is
 locally conformal parallel being the Lee form parallel,
 the corresponding four-form $\Phi$ on $S^1 \times S^7$ is a harmonic map into its sphere bundle.
\begin{proof} [Proof of Theorem \ref{cfult}]
$\,$From Proposition \ref{lcpgeneral} and Equation
\eqref{covderspin}, we have $\nabla^* \nabla \Phi =  \frac14 \|
\theta \|^2 \Phi$. Additionally, straightforward computations show
that
 \begin{eqnarray*}
\langle X \lrcorner \Phi, Y \lrcorner \Phi \rangle &= & 42 \langle
X , Y \rangle, \\
\langle X \lrcorner ( \theta^{\sharp}  \lrcorner \Phi) , Y
\lrcorner ( \theta^{\sharp}  \lrcorner \Phi) \rangle )& = & 6
\left( \| \theta \|^2  \langle X, Y \rangle  -  \theta(X) \theta
(Y)\right).
\end{eqnarray*}
Now, using these identities in Equation \eqref{scnaxy}, we obtain
$$
\langle \nabla_X \Phi, \nabla_Y \Phi \rangle = 12 \left( \| \theta
\|^2 \langle X, Y \rangle -  \theta (X) \theta (Y)\right).
$$

Finally, making use of Lemma \ref{compequi}, we get
$$
R_{(\Phi , \langle \cdot,\cdot \rangle)}= 12 \left( (d^* \theta)
\; \theta - 3 d (  \| \theta \|^2) \right).
$$
This concludes the proof.
\end{proof}
}
\end{example}

\begin{example}[Locally conformal quaternion-K\"{a}hler manifolds] {\rm
A $4n$-dimensional manifold $M$  is said to be {\it almost
quaternion-Hermitian}, if $M$ is equipped with an
$\Lie{Sp}(n)\Lie{Sp}(1)$-structure. This is equivalent to the
presence  of a Riemannian metric $\langle \cdot,\cdot \rangle$ and
a rank-three subbundle $\mathcal G$ of the endomorphism bundle
$\End TM$, such that locally $\mathcal G$ has an {\it adapted
basis} $I, J, K$ satisfying $I^2 = J^2= -1$ and $K= IJ = - JI$,
and $\langle A X, A Y \rangle = \langle X,Y \rangle $, for all
$X,Y \in T_x M$ and $A =I,J,K$. An almost quaternion-Hermitian
manifold with a global adapted basis is called an \emph{almost
hyperhermitian} manifold. In such a case the manifold is equipped
with an $\Lie{Sp}(n)$-structure.

 One may define a global,
non-degenerate four-form~$\Omega$, the {\it fundamental form}, by
the local formula
 \begin{equation} \label{fundfourform}
 \Omega = \sum_{A=I,J,K} \omega_A \wedge \omega_A,
\end{equation}
where $\omega_A (X,Y) = \langle X , A Y \rangle$, $A=I,J,K$, are
the three local K\"ahler-forms corresponding to an adapted basis.

To describe the intrinsic torsion of almost
qua\-ter\-nion-Her\-mi\-tian manifolds, it is usually used the
$\Lie{E}$-$\Lie{H}$-formalism of
\cite{Swann:symplectiques,MartinCabrera:aqh}. Thus, $\Lie{E}$ is
the fundamental representation of $\Lie{Sp}(n)$ on ${\mathbb
C}^{2n} \cong \mathbb H^n$ via left multiplication by quaternionic
matrices, considered in $\Lie{GL}(2n, {\mathbb C})$, and $\Lie{H}$
is the representation of $\Lie{Sp}(1)$ on ${\mathbb C}^2 \cong
{\mathbb H}$ given by $q . \zeta = \zeta \overline{q}$, for $q \in
\Lie{Sp}(1)$ and $\zeta \in {\mathbb H}$. An
$\Lie{Sp}(n)\Lie{Sp}(1)$-structure on a manifold $M$ gives rise to
local bundles $\Lie{E}$ and $\Lie{H}$ associated to these
representations and identifies $T M \otimes_{\mathbb R} {\mathbb
C} \cong \Lie{E} \otimes_{\mathbb C} {\Lie H}$.

The intrinsic $\Lie{Sp}(n)\Lie{Sp}(1)$-torsion $\xi$, $n>1$, is in
$ T^{\ast} M \otimes \left( {\mathfrak s}{\mathfrak p}(n) +
{\mathfrak s} {\mathfrak p}(1) \right)^{\perp}$. The decomposition
of the space of possible intrinsic torsion tensors into
irreducible $\Lie{Sp}(n)\Lie{Sp}(1)$-modules was obtained  by
Swann
 in \cite{Swann:symplectiques} and is given by
  \begin{equation} \label{swannsum}
  T^{\ast} M \otimes \left( {\mathfrak s}{\mathfrak p}(n) +
{\mathfrak s} {\mathfrak p}(1) \right)^{\perp} =  \Lambda_0^3
\Lie{E}  S^3 \Lie{H} +
  \Lie{K}  S^3 \Lie{H} + \Lie{E}  S^3 \Lie{H} +  \Lambda_0^3 \Lie{E}   \Lie{H} +
  \Lie{K}  \Lie{H} + \Lie{E}   \Lie{H},
  \end{equation}
where $\Lambda_0^3 \Lie{E}$ and  $\Lie{K}$ are certain irreducible
$\Lie{Sp}(n)$-modules and $S^3 \Lie{H}$ means the symmetric
$3$-power of $\Lie{H}$ (see \cite{MartinCabrera:aqh} for details).
If the dimension of $M$ is at least~$12$, all the modules of the
sum are non-zero. For an eight-dimensional manifold $M$, we have
$\Lambda_0^3 \Lie{E} S^3 \Lie{H}= \Lambda_0^3 \Lie{E}\Lie{H} = \{
0 \}$. Therefore, for $\dim M \geq 12 $ and $\dim M =8$, we have
respectively  $2^6 =64$ and $2^4=16$ classes of almost
quaternion-Hermitian manifolds. In \cite{MartinCabrera:aqh}, by
the identification $\xi \to - \xi \Omega = \nabla \Omega$,
explicit conditions characterising these classes were given. So
that they  are expressed in terms of $\nabla \Omega$. However,
from such conditions, it is not hard to derive  descriptions for
the $\Lie{Sp}(n)\Lie{Sp}(1)$-components of $\xi$.

A {\it locally conformal quaternion-K\"{a}hler $4n$-manifold} is an
almost quaternion Hermitian manifold such that its intrinsic
torsion $\xi$ is in $\Lie{E}   \Lie{H}$. This is equivalent to say
that
$$
\nabla_X \Omega = X^\flat \wedge (\theta^\sharp \lrcorner \Omega)
- \theta \wedge ( X \lrcorner \Omega),
$$
where the Lee form $\theta$   is given by
$$
\theta = \textstyle  \frac{1}{4(n-1)(2n+1)} \ast (\ast d \Omega
\wedge \Omega).
$$
Because of  the nature of $\Omega$, from $d \Omega = 4 \theta
\wedge \Omega$, one deduces that $\theta$ has to be closed. From
Proposition \ref{lcpgeneral} the next result follows.
\begin{theorem} \label{lckparleeqk}
For a locally  conformal quaternion-K\"{a}hler $8n$-manifold,
$\Omega^n$ is a harmonic  section of  its corresponding  sphere
bundle.
\end{theorem}

}
\end{example}

\end{document}